\mag=1000 \hsize=6.5 true in \vsize=8.7 true in
    \baselineskip=15pt
\vglue 1.5 cm


\def\move-in{\parshape=1.75true in 5true in}



    \def\qed{$\rlap{$\sqcap$}\sqcup$}


\font\tengothic=eufm10 \font\sevengothic=eufm7
\newfam\gothicfam
        \textfont\gothicfam=\tengothic
        \scriptfont\gothicfam=\sevengothic


     \font\tenmsb=msbm10              \font\sevenmsb=msbm7
\newfam\msbfam
        \textfont\msbfam=\tenmsb
        \scriptfont\msbfam=\sevenmsb
\def\Bbb#1{{\fam\msbfam #1}}


\def\PP#1{{\Bbb P}^{#1}}


\def\ref#1{[ #1]}

\def\PP#1{{\Bbb P}^{#1}}

\def\ref#1{[{\bf #1}]}

\def\a{\bigskip \par \noindent}
\def\b{\par \noindent}


\def\Proof{\noindent {\it Proof:} }

\def\th#1{\bigskip \par \noindent {\bf Theorem #1.}}






\vskip 1cm \centerline{\bf  Secant varieties to osculating varieties
of Veronese embeddings of $\PP n$.}
\bigskip

\centerline{\it A.Bernardi,
  M.V.Catalisano, A.Gimigliano, M.Id\`a}

\bigskip
\move-in {ABSTRACT: A well known theorem by Alexander-Hirschowitz
states that all the higher secant varieties of $V_{n,d}$ (the
$d$-uple embedding of $\PP n$) have the expected dimension, with few
known exceptions. We study here the same problem for $T_{n,d}$, the
tangential variety to $V_{n,d}$, and prove a conjecture, which is
the analogous of Alexander-Hirschowitz theorem, for $n\leq 9$.
Moreover. we prove that it holds for any $n,d$ if it holds for
$d=3$. Then we generalize to the case of $O_{k,n,d}$, the
$k$-osculating variety to $V_{n,d}$, proving, for $n=2$, a
conjecture that relates the defectivity of $\sigma_s(O_{k,n,d})$ to
the Hilbert function of certain sets of fat points in $\PP n$.}
\bigskip
\a {\bf  Introduction.}

\medskip
The well known Alexander-Hirschowitz  theorem (see [AH1]) states:

\th{0.1}  {\it (Alexander-Hirschowitz) Let X be a generic collection
of s 2-fat points in $\PP  n_{\kappa}$. If $(I_X)_d \subset \kappa
[x_0,...,x_n]$ is the vector space of forms of degree $d$ which are
singular at the points of $X$, then $\dim (I_X)_d = \min \{ (n+1)d,
{n+d\choose n}\}$, as expected, unless:

- $d=2, 2\leq s\leq n$;

- $n=2$, $d=4$, $s=5$;

- $n=3$, $d=4$, $s=9$;

- $n=4$, $d=3$, $s=7$;

- $n=4$, $d=4$, $s=14$.}

\bigskip
Notice that with
``$m$-fat point at $P\in \PP n$" we mean the scheme defined by the ideal
$I_P^m\subset \kappa
[x_0,...,x_n]$.

 An equivalent
reformulation of the theorem is in the language of higher secant
varieties; let $V_{n,d}\subset \PP N$, with $N={n+d\choose n}-1$, be
the $d$-ple (Veronese) embedding  of $\PP n$, and let $\sigma_s
(V_{n,d})$ be its $(s-1)^{th}$ higher secant variety, that is, the
closure of the union of the $\PP {s-1}$'s which are $s$-secant to
$V_{n,d}$. Then Theorem 0.1 is equivalent to:

\medskip

\th{0.2}  {\it All the higher secant varieties $\sigma_s (V_{n,d})$
have the expected dimension $\min \{ s(n+1)-1, {n+d\choose n}-1\}$,
unless $s,n,d$ are as in the exceptions of Theorem 0.1.}

An application of the theorem is in terms of the Waring problem for
forms (or of the decomposition of a supersymmetric tensor), in fact
Theorem 0.1 gives that the general form of degree $d$ in $n+1$
variables can be written as the sum of $\lceil {1\over
n+1}{n+d\choose d}\rceil$ $d$th powers of linear forms, with the
same list of exceptions (e.g. see [Ge] or [IK]).

\medskip

In [CGG] a similar problem has been studied, namely whether the
dimension of $\sigma_s (T_{n,d})$ is the expected one or not, where
$T_{n,d}$ is the tangential variety of the Veronese variety
$V_{n,d}$. This too translates into a problem of representation of
forms: the generic form parameterized by $\sigma_s (T_{n,d})$ is a
form $F$ of degree $d$ which can be written as $F = L_1^{d-1}M_1 +
... + L_s^{d-1}M_s$, where the $L_i,M_i$'s are linear forms.

 The following conjecture was stated in [CGG]:

\medskip
{\bf Conjecture 1:} {\it The secant variety} $\sigma_{s}(T_{n,d})$
{\it has the expected dimension, $\min \{ 2sn+s-1, {n+d\choose
n}-1\}$, unless:}
\par
$i)\;\;$  $d=2$ {\it ,}\ $2\leq 2s < n${\it ;}
\par
$ii)\;$ $d=3${\it ,} $s=n=2,3,4${\it .}
\par
\medskip

In the same paper the conjecture was proved for $d=2$ (any $s,n$)
and for $s\leq 5$ (any $d,n$), while in [B] it is proved for $n=2,3$
(any $s,d$).

\medskip In [CGG](via inverse systems) it is shown that $\sigma_s (T_{n,d})$
is defective if and only if a certain 0-dimensional scheme
$Y\subset \PP n$ does not impose independent conditions to forms
of degree $d$ in $R:= \kappa[x_0,...,x_n]$. The scheme $Y= Z_1\cup
...\cup Z_s$ is supported at $s$ generic points $P_1,...,P_s \in
\PP n$, and at each of them the scheme $Z_i$ lies between the
2-fat point and the $3$-fat point on $P_i$ (we will call $Z_i$ a
$(2,3,n)$-scheme, for details see section 1 below).

Hence Conjecture 1 can be reformulated in term of $(I_Y)_d$ having
the expected dimension, with the same exceptions, in analogy with
the statement of Theorem 0.1.

Theorem 0.1 has been proved thanks to the Horace differential Lemma
({AH2}, Proposition 9.1; see also here Proposition 1.5) and an
induction procedure which has a delicate beginning step for $d=3$;
different proofs for this case are in [Ch1], [Ch2] and in the more
recent [BO], where an excellent history of the question can be
found.

Also the proof of Conjecture 1 presents the case of $d=3$ as a
crucial one; the first main result in this paper (Corollary
2.5) is to prove that if Conjecture 1 holds for $d=3$, then it holds
also for $d\geq 4$ (and any $n,s$).  The procedure we use is based
on Horace differential Lemma too.

We also prove Conjecture 1 for all $n\leq 9$, since
with that hypothesis we can check the case $d=3$ by making use of
{COCOA} (see Corollary 2.4).

\bigskip

A more general problem can be considered (see also [BCGI]): let
$O_{k,n,d}$ be the $k$-osculating variety to $V_{n,d}\subset \PP N$,
and  study its $(s-1)^{th}$ higher secant variety $\sigma_s
(O_{k,n,d})$. Again, we are interested in the problem of
determining all $s$ for which $\sigma_s (O_{k,n,d})$ is defective,
i.e. for which its dimension is strictly less than its expected
dimension (for precise definitions and setting of the problem, see
Section 1 of the present paper and in particular Question
Q(k,n,d)).

\medskip Also in this general case we found in [BCGI] (via inverse systems)
 that $\sigma_s (O_{k,n,d})$ is defective if and only if a certain 0-dimensional scheme
$Y\subset \PP n$ does not impose independent conditions to forms of
degree $d$ in $R:= \kappa[x_0,...,x_n]$. The scheme $Y= Z_1\cup
...\cup Z_s$ is supported at $s$ generic points $P_1,...,P_s \in \PP
n$, and at each of them the ideal of the scheme $Z_i$ is such that
$I_{P_i}^{k+2} \subset I_{Z_i} \subset I_{P_i}^{k+1}$ (for details
see Lemma 1.2 below).

The following  (quite immediate) lemma ([BCGI] 3.1) describes what
can be deduced about the postulation of the scheme $Y$ from
information on fat points:

\a {\bf Lemma 0.3.} {\it Let $P_1,...,P_s$ be generic points in $\PP
n$, and set $X:= (k+1)P_1\cup ...\cup (k+1)P_s$, $T:= (k+2)P_1\cup
...\cup (k+2)P_s$. Now let $Z_i$ be a 0-dimensional scheme supported
at $P_i\,$, $(k+1)P_i\subset Z_i\subset (k+2)P_i$, and set $Y:=
Z_1\cup ...\cup Z_s$. Then, $Y$ is regular in degree $d$ if
$h^1({\cal I}_{T}(d))=0$ or if $h^0({\cal I}_{X}(d))=0$. \b
Moreover, $Y$ is not regular in degree $d$ if \b (i) $h^1({\cal
I}_{X}(d))> max \{0, \deg(Y) - {d+n \choose n} \}$,
 \b or if \b (ii)
$\;h^0({\cal I}_ {T}(d))>max \{0, {d+n \choose n}-\deg(Y) \}$.}

\a All cases studied in [BCGI] lead us to state the following:

\medskip \noindent {\bf Conjecture 2a.} {\it The secant variety $\sigma_s (O_{k,n,d})$ is
defective if and  only if $Y$ is as in case (i) or (ii) of the Lemma
above.}
\medskip
The conjecture amounts to saying that $I_Y$ does not have the
expected Hilbert function in degree $d$ only when ``forced" by the
Hilbert function of one of the fat point schemes $X$, $T$.

Notice that (i), respectively (ii), obviously implies that
 $X$, respectively $T$, is defective. Hence, if Conjecture 2a holds
and  $Y$ is defective in degree $d$, then  either $T$ or $X$ are
defective in degree $d$ too, and the defectivity of $Y$ is either
given by the defectivity of $X$ or forced by the high defectivity of
$T$.

Thus if the conjecture holds, we have another occurrence of the
"ubiquity" of fat points: the problem of $\sigma_s (O_{k,n,d})$
having the right dimension reduces to a problem of computing the
Hilbert function in degree $d$ of two schemes of $s$ generic fat
points in $\PP n$, all of them having multiplicity $k+1$,
respectively $k+2$.

\medskip
 In [BC]and [BF] the conjecture is proved in $\PP 2$ for $s\leq 9$.

\medskip
 Notice that the Conjecture 2a implies the following one, more
geometric, which relates the defectivity of $\sigma_s (O_{k,n,d})$
to the dimensions of the $k^{th}$ and the $(k+1)^{th}$ osculating
space at a generic point of the $(s-1)^{th}$ higher secant variety
of the Veronese variety $\sigma_s (V_{n,d})$:

\medskip

 \noindent {\bf Conjecture 2b.} {\it If the secant variety $\sigma_s (O_{k,n,d})$ is
defective then at a generic point $P\in \sigma_s (V_{n,d})$, either the $k^{th}$ osculating
 space $O_{k,\sigma_s (V_{n,d}),P}$ does not have dimension
$\min \{s{k+n \choose n}-1,{d+n \choose n}-1\}$, or the $(k+1)^{th}$ osculating
space $O_{k+1,\sigma_s (V_{n,d}),P}$  does not have dimension
$\min \{s{k+n+1 \choose n}-1, {d+n \choose n}-1\}.$}
\medskip \noindent
 The implication follows from the fact
 that (see  [BBCF]) for $P\in <P_1,...,P_s>$:
 $$O_{k,\sigma_s (V_{n,d}),P}=
 <O_{k,V_{n,d},P_1},\;O_{k,V_{n,d},P_2},\;\dots,O_{k,V_{n,d},P_s}\;>.$$

The other main result in this paper is Theorem 3.5, which proves
Conjecture 2a for $n=2$.

 \bigskip
\a {\bf  Section 1: Preliminaries and Notations.}

\a In this paper we will always work over a field $\kappa$ such
that $\kappa=\overline \kappa$ and char$\kappa = 0$.

\a{\bf 1.1 Notations.}
\par  (i) If  $P\in  \PP n$ is
a point and $I_P$ is the ideal of $P$ in  $\PP n$, we denote by $mP$
the fat point of multiplicity $m$ supported at $P$, i.e. the scheme
defined by the ideal $I_P^m$.

\par (ii) Let $X\subseteq \PP N$ be a closed irreducible
projective variety;  the $(s-1)^{th}$ {\it higher secant variety} of
$X$ is the closure of the union of all linear spaces spanned by $s$
points of $X$, and it will be denoted by $\sigma_s(X)$.

\par  \medskip (iii) Let $X \subset \PP N$ be a variety, and let $P\in X$ be
a smooth point; we define {\it the $k^{th}$ osculating space to X at
P} as the linear space generated by $(k+1)P\cap X$ (i.e. by the
$k^{th}$ infinitesimal neighbourhood of $P$ in $X$) and we denote it
by $O_{k,X,P}$; hence $O_{0,X,P}=\{P\}$, and $O_{1,X,P}=T_{X,P}$,
the projectivised tangent space to $X$ at $P$.
\smallskip \b
Let $U\subset X$ be the dense set of the smooth points where
$O_{k,X,P}$ has maximal dimension. The {\it $k^{th}$ osculating
variety to X } is defined as:
$$ O_{k,X}= \overline {\bigcup _{P\in U}O_{k,X,P}}. $$

\par  \medskip (iv) We denote by $V_{n,d}$ the $d$-uple Veronese embedding of $\PP n$,
i.e. the image of the map  defined by the linear system of all forms
of degree $d$ on $\PP n$: $\nu_d: \PP n \rightarrow \PP N$, where
$N={n+d\choose n}-1$.

\par  \medskip (v)  We denote the $k^{th}$ osculating
variety to the Veronese  variety  by $O_{k,n,d}:=O_{k,V_{n,d}}$.
When $k=1$, the osculating variety is called {\it tangential
variety} and it is denoted by $T_{n,d}$.
\par Hence, the
$(s-1)^ {th}$ higher secant variety of the $k^{th}$ osculating
variety to the Veronese  variety $V_{n,d}$ will be denoted by
$\sigma _s(O_{k,n,d})$.

\a Since the case $d \leq k$ is trivial, and the description for $k = 1$ given in [CGG],
 together with [BCGI, Proposition 4.4] describe the case $d = k+1$ completely,
 from now on we make the general assumption, which will be implicit in the rest of the  paper, that $d \geq k+2$.

\a It is easy to see ([BCGI] 2.3) that the dimension of $O_{k,n,d}$ is always
the expected one, that is, ${\dim} O_{k,n,d} = {\rm min} \{N,\;
n+{k+n\choose n}-1\}$. The expected dimension for $\sigma
_s(O_{k,n,d})$ is:
$$ {\rm expdim}\, \sigma _s(O_{k,n,d}) = \min \{\, N,\, s(n+{k+n\choose n}-1)+s-1\}$$

\noindent (there are $\infty ^{s({\rm dim}O_{k,n,d})}$ choices of
$s$ points on $O_{k,n,d}$, plus $\infty ^{s-1}$ choices of a point
on the $\PP {s-1}$ spanned by the $s$ points; when this number is
too big,  we expect that $\sigma _s(O_{k,n,d}) = \PP N$). \b When $
\dim \sigma _s(O_{k,n,d}) <{\rm expdim}\, \sigma _s(O_{k,n,d})$, the
osculating variety is said to be {\it defective}.

\par  \medskip In [BCGI], taking into account that the cases with $n=1$ can be easily
described, while if $n\geq 2$ and $d=k$ one has ${\rm dim}\, \sigma
_s(O_{k,n,d})=N$, we raised the following question:

\a {\bf Question Q(k,n,d)}: {\it For all $k,n,d$ such that $d \geq
k+1$, $n\geq 2$, describe all $s$ for which $\sigma_s (O_{k,n,d})$
is defective, i.e.}
$$\dim \sigma_s (O_{k,n,d}) < \min \{\, N,\, s(n+{k+n\choose n}-1)+s-1\}=\min \{\,
{d+n\choose n}-1,\, s{k+n\choose n}+sn-1\}.$$

We were able to answer the question for $s,n,d,k$ in several ranges,
thanks to the following lemma (see [BCGI] 2.11 and results of
Section 2):

\a {\bf Lemma {1.2}} {\it For any $k,n,d \in \Bbb N$ such that
$n\geq 2$, $d\geq k+1$, there exists a 0-dimensional subscheme
$Z=Z(k,n)\in \PP n$ depending only from $k$ and $n$ and not from
$d$,  such that:

\par (a) $Z$ is supported on a point $P$, and one has:
$$(k+1)P \subset Z(k,n) \subset (k+2)P, \qquad with  \qquad l(Z)={k+n \choose n}+n;$$

(b) denoting by $Y=Y(k,n,s)$ the generic union in $\PP n$ of
$Z_1,\dots,Z_s$ where $Z_i\cong Z$ for $i=1,\dots,s$, then
$$\dim \sigma_s (O_{k,n,d}) = {\rm expdim} \,\sigma_s (O_{k,n,d})-h^0({\cal I}_{Y}(d))+
\max \{0,{d+n\choose n}-l(Y)\}$$

\noindent In particular, $\sigma_s (O_{k,n,d})$ is not defective if
and only if $Y$ is regular in degree $d$, i.e. $h^0({\cal
I}_Y(d))\cdot h^1({\cal I}_Y(d))=0$.}

\a The homogeneous ideal of this 0-dimensional scheme $Z$ is defined
in [BCGI] 2.5 through inverse systems, so we don't have an explicit
geometric description of it in the general case. Anyway, for $k=1$
it is possible to describe it geometrically as follows (see [CGG]
Section 2):

\a {\bf  Definition 1.3.} Let $P$ be a point in $\PP n$, and $L$ a
line through $P$; we say that a 0-dimensional scheme $X\subset \PP
n$ is a $(2,3,n)$-scheme supported on $P$ with direction $L$ if
$I_X=I_P^3+I_L^2$. Hence, the length of a $(2,3,n)$-point is $2n+1$.
The scheme $Z(1,n)$ of Lemma 1.2 is a $(2,3,n)$-scheme.

\par \smallskip We say that a
subscheme of $\PP n$ is a generic union of $s$ $(2,3,n)$-schemes if
it is the union of $X_1,\dots,X_s$  where $X_i$ is a
$(2,3,n)$-scheme supported on $P_i$ with direction $L_i$, with
$P_1,\dots,P_s$ generic points and $L_1,\dots,L_s$ generic lines
through $P_1,\dots,P_s$.

\par \bigskip We are going to use these schemes in Section 2, so we need to know more
about them; but first we recall the Differential Horace Lemma of
\ref {AH2}, writing it in the context where we shall use it.

\a {\bf  Definition 1.4.} In the algebra of formal functions $\kappa
[[{\bf x},y]]$, where ${\bf x} = (x_1,...,x_{n-1})$, a {\it
vertically graded} (with respect to $y$) ideal is an ideal of the
form:
$$
I = I_0 \oplus I_1y \oplus ... \oplus I_{m-1}y^{m-1}\oplus (y^m)
$$
where for $i = 0,...,m-1$, $I_i\subset \kappa [[{\bf x}]]$ is an
ideal.
\par
Let $Q$ be a smooth $n$-dimensional integral scheme, let $K$ be a
smooth irreducible divisor on $Q$. We say that $Z \subset Q$ is a
{\it vertically graded subscheme} of $Q$ with base $K$ and support
$z\in K$, if $Z$ is a 0-dimensional scheme with support at the point
$z$ such that there is a regular system of parameters $({\bf x},y)$
at $z$ such that $y=0$ is a local equation for $K$ and the ideal of
$Z$ in $\widehat {\cal O}_{Q,z} \cong \kappa [[{\bf x},y]]$ is
vertically graded.
\par
\medskip

\par \noindent Let $Z\i Q$ be a vertically graded subscheme with base $K$, and $p\geq 0$
be a fixed integer; we denote by  $Res^{p}_K(Z)\i Q$ and
$Tr^{p}_K(Z)\i K$  the closed subschemes defined, respectively, by
the ideals:

$$ {\cal I}_{Res^{p}_K(Z)} := {\cal I}_{Z}+
({\cal I}_{Z}:{\cal I}^{p+1}_K){\cal I}^{p}_K, \qquad \qquad \qquad
{\cal I}_{Tr^{p}_K(Z),K} := ({\cal I}_{Z}:{\cal I}^{p}_K)\otimes
{\cal O}_K. $$

\par \noindent In $Res^{p}_K(Z)$ we take away from $Z$ the $(p+1)^{th}$ ``slice"; in
$Tr^{p}_K(Z)$ we consider only the $(p+1)^{th}$ ``slice". Notice
that for $p=0$ we get the usual trace and residual schemes:
$Tr_K(Z)$ and  $Res_K(Z)$.

\par Finally, let $Z_1,...,Z_r\i Q$ be vertically graded subschemes with base $K$
and support $z_i$, $Z=Z_1\cup...\cup Z_r$, and ${\bf
p}=(p_1,...,p_r)\in {\Bbb N}^r$.
\par We set:
$$Tr^{\bf p}_K(Z):= Tr^{p_1}_K(Z_1)\cup ... \cup Tr^{p_r}_K(Z_r),
\quad Res^{\bf p}_K(Z):= Res^{p_1}_K(Z_1)\cup ... \cup
Res^{p_r}_K(Z_r).$$

\b {\bf  Proposition 1.5.} (Horace differential Lemma, \ref {AH2}
Proposition 9.1) {\it Let $H$ be a hyperplane in $\PP n$ and let
$W\i \PP n$ be a 0-dimensional closed subscheme.
\par
Let $S_1,...,S_r,\, Z_1,...,Z_r$ be $0-$dimensional irreducible
subschemes of $\PP n$ such that $S_i\cong Z_i$, $i=1,...,r$, $Z_i$
has support on $H$ and is vertically graded with base $H$, and the
supports of $S=S_1\cup ....\cup S_r$ and $Z=Z_1\cup ...\cup Z_r$ are
generic in their respective Hilbert schemes. Let ${\bf
p}=(p_1,...,p_r)\in {\Bbb N}^r$. Assume:
\par \noindent a)  $H^0({\cal I}_{Tr_HW\cup Tr^{\bf p}_H(Z),H}(n))=0$ and

\par \noindent  b) $H^0({\cal I}_{Res_HW\cup Res^{\bf p}_H(Z)}(n-1))=0$,
\par
\noindent then
$$ H^0({\cal I}_{W\cup S}\,(n))=0. $$ }

 \a {\bf Definition 1.6.} A {\it $2$-$jet$} is  a 0-dimensional scheme
$J\subset \PP n$ with support at a point $P\in \PP n$ and degree 2;
namely the ideal of $J$ is of type: $I_P ^2 + I_L$, where $L\subset
\PP n$ is a line containing $P$.  We will say that $J_1,...,J_s$ are
$generic$ in $\PP n$, if the points $P_1,...,P_s$ are generic in
$\PP n$ and $L_1,\dots,L_s$ are generic lines through
$P_1,\dots,P_s$.

\a {\bf Remark 1.7.} Let $X\subset \PP n$ be a $(2,3,n)$-scheme
supported at $P$ with direction $L$ and $(y_1,\dots,y_n)$ be local
coordinates around $P$, such that $L$ becomes the $y_n$-axis; then,
$I_X=(y_1y_n^2,\dots,y_{n-1}y_n^2,
y_n^3,y_1^2,y_1y_2,\dots,y_{n-1}^2)$ ($y_n$ appears only in the
first $n$ generators). Let $H$, respectively $K$, be a hyperplane
through $L$, respectively transversal to $L$; then, we can assume
$I_H=(y_{n-1})$, respectively $I_K=(y_{n})$. We now compute
$Res^{p}_H(X)$ and $Tr^{p}_H(X)$. One has:

\medskip \par a) $Res _H X= Res^{0}_H(X)$, $I_{Res_H(X)} =
(I_X:y_{n-1})=(y_1,\dots,y_{n-1},y_n^2)$, hence $Res _H X$ is a
2-jet lying on $L$;
\medskip \par b) $Tr_H(X)=Tr^0_H(X)$, $I_{Tr_H(X)} =
I_X+(y_{n-1})=(y_1y_n^2,\dots,y_{n-2}y_n^2,
y_n^3,y_1^2,y_1y_2,\dots,y_{n-2}^2)$, hence $Tr_H(X)$ is a
$(2,3,n-1)$-scheme of $H$.

\medskip \par  Hence the scheme $X$ as a
vertically graded scheme with base $H$ has only two layers (strata);
in other words, $Tr^p_H(X)$ is empty for $p>1$, and $Res^1_H(X)$ is
a (2,3,n-1)-scheme of $H$, while $Tr^1_H(X)$ is a 2-jet lying on
$L$.

\medskip \par   Now we want to compute $Res^{p}_K(X)$ and
$Tr^{p}_K(X)$. Consider first:

\medskip \par b) $I_{Tr_K(X)} = I_X+(y_n)=(y_n,y_1^2,y_1y_2,\dots,y_{n-1}^2)$,
hence $Tr_H(X)$ is a 2-fat point of $K\cong \PP {n-1}$,

\medskip \par a)  $I_{Res_KX} =
(I_X:y_n)=(y_1y_n,\dots,y_{n-1}y_n,y_n^2,y_1^2,y_1y_2,\dots,y_{n-1}^2)$,
hence $Res_KX$ is a 2-fat point of $\PP n$.

\medskip \par  So the scheme $X$, as a
vertically graded scheme with base $K$, has only three layers
(strata); the 0-layer is $Tr_K(X)=Tr^0_K(X)$, the 1-layer is the
0-layer of $Res _K X= Res^{0}_K(X)$, hence it is again a 2-fat point
of $K\cong \PP {n-1}$, and the 2-layer is  the 1-layer of $Res _K
X$, hence it is a point of $\PP n$. In other words, $Tr^p_H(X)$ is
empty for $p>2$, $Res^1_K(X)$ is a a 2-fat point of $\PP n$, while
$Tr^1_K(X)$ is a 2-fat point of $K$; $Res^2_K(X)$ is a  2-fat point
of $K$  doubled in a direction transversal to $K$ (i.e.,
$I_{Res^2_K(X)}=(y_n^2,y_1^2,y_1y_2,\dots,y_{n-1}^2)$), while
$Tr^2_K(X)$ is a point of $\PP n$.

\bigskip We will use in the sequel the fact that by adding $s$ generic 2-jets to any 0-dimensional
scheme $Z\subset \PP n$ we impose a maximal number of independent
conditions to forms in $I_Z(d)$, for all $d$. This is probably
classically known, but we write a proof here for lack of a
reference:

\a {\bf Lemma {1.8}} {\it Let $Z\subseteq \PP n$ be a scheme, and
let $J\subset \PP n$ be a generic 2-jet. Then:
$$h^0({\cal I}_{Z\cup J}(d)) = \max \{h^0({\cal I}_{Z}(d))-2,0\}.$$}

\b \Proof  Let $P$ be the support of $J$; then we know that
$h^0({\cal I}_{Z\cup P}(d)) = \max \{h^0({\cal I}_Z(d))-1,0\}$, so
if $h^0({\cal I}_Z(d))\leq 1$ there is nothing to prove. Let
$h^0({\cal I}_Z(d)) \geq 2$, then $h^0({\cal I}_{Z\cup P}(d)) =
h^0({\cal I}_Z(d))-1 \geq 1$. Since $J$ is generic, if $ h^0({\cal
I}_{Z\cup J}(d)) = h^0({\cal I}_{Z\cup P}(d))$, then every form of
degree $d$ containing $Z\cup P$ should have double intersection with
almost every line containing $P$, hence it should be singular at
$P$. This means that when we force a form in the linear system
$\vert H^0({\cal I}_{Z}(d))\vert$ to vanish at $P$, then we are
automatically imposing to the form to be singular at $P$, and this
holds for $P$ in a dense open set of $\PP n$, say $U$. If the form
$f$ is generic in $\vert H^0({\cal I}_{Z}(d))\vert$, its zero set
$V$ meets $U$ in a non empty subset of $V$, so $f$ is singular at
whatever point $P'$ we choose in $V\cap U$, and this means that the
hypersurface $V$ is not reduced. Since the dimension of the linear
system $\vert H^0({\cal I}_{Z}(d))\vert$ is at least 2, this is
impossible by Bertini Theorem (e.g. see \ref {J}, Theorem 6.3).

\hfill\qed

\bigskip
Let $Z \subseteq \PP n$ be a zero-dimensional scheme; the following
simple Lemma gives a criterion for adding to $Z$ a scheme $D$ which
lies on a smooth hypersurface ${\cal F} \subseteq \PP n$ and is made
of $s$ generic 2-jets on ${\cal F}$, in such a way that $D$ imposes
independent conditions to forms of a given degree in the ideal of
$Z$ (see Lemma 4 in \ref {Ch1} and Lemma 1.9 in \ref {CGG2} for the
case of simple points on a hypersurface).

\a {\bf Lemma {1.9}} {\it Let $Z \subseteq \PP n$ be a zero
dimensional scheme.  Let ${\cal F} \subseteq \PP n$ be a smooth
hypersurface of degree $d$ and let $Z'=Res_{\cal F} Z$. Let
$P_1,\dots,P_s$ be generic points on $\cal F$, let $L_1,\dots,L_s$
lines with $P_i\in L_i$, and such that each line $L_i$ is generic in
$T_{P_i}({\cal F})$; let $J_i$ be the 2-jet with support at $P_i$
and contained in $L_i$. We denote by $D_s=J_1\cup \dots \cup J_s$
the union of these $s$ 2-jets  generic in ${\cal F}$.
\par i) If $\dim (I_{Z+D_{s-1}})_t \geq \dim (I_{Z'})_{t-d} + 2$,
then $ \dim (I_{Z+D_s})_t =  \dim (I_{Z})_t -2s;$}
\par
{\it ii) if $ \dim (I_{Z'})_{t-d}=0$ and $\dim (I_{Z})_t\leq 2s$,
then $\dim (I_{Z+D_s})_t = 0$.}

\medskip \Proof \  {\it i)} By induction on $s$. If $s=1$, by assumption $\dim (I_{Z})_t
\geq \dim (I_{Z'})_{t-d}+2$,  hence in the exact sequence $0\to
H^0({\cal I}_{Z'}(t-d))\buildrel {\phi}\over \to H^0({\cal
I}_Z(t-d))\to H^0({\cal I}_{Z\cap{\cal F},{\cal F}}(t))\to\dots$ the
cokernel of the map $\phi$ has dimension at least 2 and so
$(I_{Z})_t$ cuts on ${\cal F}$ a linear system (i.e. $\vert
H^0({\cal I}_{Z\cap{\cal F},{\cal F}}(t)) \vert$) of (projective)
dimension $\geq 1$. We have $\dim (I_{Z+P_1})_t =\dim (I_{Z})_t -1$,
since otherwise each hypersurface in $\vert (I_{Z})_t \vert$ would
contain the generic point $P_1$ of ${\cal F}$, that is, would
contain ${\cal F}$.
\par Assume $\dim (I_{Z+J_1})_t =\dim (I_{Z+P_1})_t =\dim (I_{Z})_t -1$; this means that
if we impose to $S\in \vert (I_{Z})_t \vert$  the passage through
$P_1$ automatically we impose to $S$ to be tangent to $L_1$ at
$P_1$, and $L_1$ being generic in $T_{P_1}({\cal F})$, this means
that each $S$ passing through $P_1$ is tangent to ${\cal F}$ at
$P_1$. Let's say that this holds for $P_1$ in the open not empty
subset $U$ of ${\cal F}$; for $S$ generic in $\vert (I_{Z})_t
\vert$, $U'=S\cap {\cal F} \cap U$ is not empty, hence the generic
$S$ is tangent to ${\cal F}$ at each $P\in U'$. This means that
$\vert (I_{Z})_t \vert$ cuts on ${\cal F}$ a linear system of
positive dimension whose generic element is generically non reduced,
and this is impossible, by Bertini Theorem (e.g. see \ref {J},
Theorem 6.3).
\par \smallskip
Now let $s>1$.  Since $\dim (I_{Z+D_{s-2}})_t\geq \dim
(I_{Z+D_{s-1}})_t >  \dim (I_{Z'})_{t-d}$ by assumption, and
$Res_{\cal F} (Z+D_{s-1})=Z'$, the case $s=1$ gives $\dim (I_{Z+D_s}
)_t = \dim(I_{Z+D_{s-1}} )_t -2$.
 So, by the induction hypothesis, we get
$$
\dim (I_{Z+D_s} )_t = (\dim (I_{Z})_t -2(s-1))-2 =\dim (I_{Z})_t
-2s.
$$
{\it ii)}\ Assume first $\dim (I_{Z})_t \leq 2$; it is enough to
prove $\dim (I_{Z+J_1})_t = 0$ since then also $\dim (I_{Z+D_s})_t =
0$. If $\dim (I_{Z})_t = 2$ this follows by $i)$ and if $\dim
(I_{Z})_t = 0$ this is trivial. If $\dim (I_{Z})_t = 1$, then if
$\dim (I_{Z+P_1})_t = 0$ we are done. If $\dim (I_{Z+P_1})_t = 1$,
then by the genericity of $P_1$ we have that the unique $S$ in the
system contains ${\cal F}$, i.e. $S={\cal F}\cup G$, but then
$Z'\subseteq G$, which contradicts $\dim (I_{Z'})_{t-d}=0$.

Otherwise, let $\dim (I_{Z})_t = 2v+\delta \geq 3$, $\delta =0,1$.
If $\delta =0$, then $\dim (I_{Z+D_{v-1}})_t \geq 2 = \dim
(I_{Z'})_{t-d}+2$, and by {\it i)} we get  $\dim (I_{Z+D_{v}})_t =
\dim (I_{Z})_t -2v = 0$, and, since $s \geq v$, it follows that
$\dim (I_{Z+D_{s}})_t =0$.

If $\delta =1$, then  $\dim (I_{Z+D_{v-1}})_t \geq 3 \geq \dim
(I_{Z'})_{t-d}+2$, and, by $i)$, $\dim (I_{Z+D_{v-1}})_t = 3$ and
$\dim (I_{Z+D_{v}})_t = \dim (I_Z)_t - 2v = 1$.  Notice that the
only element in $(I_{Z+D_{v}})_t$ cannot have ${\cal F}$ as a fixed
component, otherwise we would have $\dim (I_{Z'})_{t-d}=1$ and not
$=0$; hence $\dim (I_{Z+D_v+P_{v+1}})_t = 0$ and so, since $2s\geq
2v+1$ and $D_v \cup P_{v+1}\subset D_s$, $\dim (I_{D_s})_t = 0$.

\hfill\qed

\a Now we give a Lemma which will be of use in the proof of Theorem
2.2.

\bigskip

{\bf Lemma {1.10}} {\it Let $R \subseteq \PP n$ be a zero
dimensional scheme contained in a $(2,3,n)$-scheme with $r=\deg Y
\leq 2n$; assume moreover that, if $r\geq n+1$, then $R$ is a flat
limit of the union of a 2-fat point of $\PP n$ and of a scheme
(eventually empty) contained in a 2-fat point of a $\PP {n-1}$, and
that, if $r\leq n$, then $R$ is contained in a 2-fat point of a $\PP
{n-1}$. Then, there exists a flat family for which R is a special
fiber and the generic fiber is the generic union in $\PP n$ of
$\delta$ 2-fat points, $h$ 2-jets and $\epsilon$ simple points,
where $r=(n+1)\delta + 2h+ \epsilon$, $0\leq \delta \leq 1$, $0\leq
\epsilon \leq 1$, and $2h+ \epsilon\leq n$. }

\medskip \Proof \  In the following we denote by $2_tP$ a 2-fat point of a linear variety
$K\subseteq \PP n$, $K\cong \PP t$. We first notice that if $A$ is a
subscheme of $2_nP$ with $\deg A= n$ then $A$ is a scheme of type
$2_{n-1}P$. The proof is by induction on $n$: if $n=2$, the
statement is trivial since the only scheme of degree 2 in $\PP 2$ is
a 2-jet, i.e. a $2_1P$. Now assume the assertion true for $n-1$, let
$A$ be a subscheme of $2_nP$ with $\deg A= n$ and let $H$ be a
hyperplane through the support of $A$. Since $\deg 2_nP\cap H=n$, we
have $n-1\leq \deg A\cap H\leq n$. If $\deg A\cap H= n$ then
$A=2_{n-1}P$ and we are done. If $\deg A\cap H= n-1$ then $Res_H A$
is a simple point, and by induction $A\cap H=2_{n-2}P$. Hence there
is a hyperplane $K$ such that $A\cap H$ is a 2-fat point of $H\cap
K$, and working for example in affine coordinates, it is easy to see
that $A$ is a 2-fat point of the $\PP {n-1}$ generated by $H\cap K$
and a normal direction to $H$.
\medskip \par In order to prove the Lemma, it is enough to prove that the generic union
in $\PP n$ of $h$ 2-jets and $\epsilon$ simple points, with $0\leq
\epsilon \leq 1$ and $2h+ \epsilon\leq n$, specializes to any
possible subscheme $M$ of a scheme of type $2_{n-1}P$: in fact, if
$r\leq n$ we are done, if $r\geq n+1$, the collision of a $2_nP$
with $M$ gives $R$.
\par By induction on $n$: if $n=2$, the statement is trivial. Let us now consider the
generic union of $h$ 2-jets and $\epsilon$ simple points in $\PP n$,
with $0\leq \epsilon \leq 1$ and $2h+ \epsilon\leq n$. We have two
cases.
\par Case 1: if $2h+ \epsilon\leq n-1$, we specialize everything inside a hyperplane $H$
where, by induction assumption, this scheme specializes to any
possible subscheme of a scheme of type $2_{n-2}P$, i.e., to  any
possible subscheme of degree $\leq n-1$ of a scheme of type
$2_{n-1}P$.
\par Case 2: If $2h+ \epsilon = n$, we have to show that the generic union of $h$
2-jets and $\epsilon$ simple points specializes to a scheme
$2_{n-1}P$.
\par If $n$ is odd, then $h={n-1 \over 2}$ and $\epsilon =1$; by induction assumption,
${n-1 \over 2}$ 2-jets specialize to a scheme of type $2_{n-2}P$,
and the generic union of the last one with a simple point
specializes to a scheme of type $2_{n-1}P$.
\par If $n$ is even, then $h={n \over 2}$ and $\epsilon =0$;  by induction assumption,
${n \over 2}-1$ 2-jets specialize to a scheme of degree $n-2$
contained in a scheme of type $2_{n-2}P$, which is a $2_{n-3}P$, so
it is enough to prove that the generic union of the last one with a
2-jet specializes to a scheme of type $2_{n-1}P$.
\par In affine coordinates $x_1,\dots,x_n$, let $x_{n-2}=x_{n-1}=x_{n}=0$ be the linear
subspace containing $2_{n-3}P$, so that
$I_{2_{n-3}P}=(x_1,\dots,x_{n-3})^2 \cap (x_{n-2},x_{n-1},x_{n})$,
and let $(x_1,\dots,x_{n-3}, x_{n-2}-a,x_{n-1}^2,x_{n})$ be the
ideal of a 2-jet moving along the $x_{n-2}$-axis; then it is
immediate to see that the limit for $a\rightarrow 0$ of
$(x_1,\dots,x_{n-3})^2 \cap (x_{n-2},x_{n-1},x_{n})\cap
(x_1,\dots,x_{n-3}, x_{n-2}-a,x_{n-1}^2,x_{n})$ is
$(x_1,\dots,x_{n-1})^2 \cap (x_{n})$, which is the ideal of a
$2_{n-1}P$.

\hfill\qed

\bigskip \noindent {\bf 2. On Conjecture 1.}

\medskip

We want to  study $\sigma_{s}(T_{n,d})$, and we have seen that its
dimension is given by the Hilbert function of $s$ generic
$(2,3,n)$-points in $\PP n$.

\a {\bf Definition 2.0} For each $n$ and $d$ we define  $s_{n,d}$,
$r_{n,d} \in \Bbb{N}$ as the two positive integers such that
$${d+n\choose
n}=(2n+1)s_{n,d}+r_{n,d}, \qquad 0\leq r_{n,d} <2n+1.$$

In the following we denote by $X_{s,n}\subset \PP n$ the zero
dimensional scheme union of $s$ generic $(2,3,n)$-schemes
$A_1,\dots,A_s$. We also denote by $X_{s_{n,d}}$ the scheme
$X_{s,n}$, with $s = s_{n,d}$.  Hence $X_{s_{n,d}}$ is the union of
the maximum number of generic $(2,3,n)$-points that we expect to
impose independent conditions to forms od degree $d$. We will also
use $X_{s_{n,d}+1}$ to indicate $X_{s+1,n}$ when $s = s_{n,d}$.

\par  With $Y_{n,d}\subset \PP n$ we denote a scheme generic union of $X_{s_{n,d}}$ and
$R_{n,d}$, where $R_{n,d}$ is a zero dimensional scheme contained in
a $(2,3,n)$-point, with $\deg(R_{n,d})= r_{n,d}$.
\medskip
 A 0-dimensional subscheme $A$ of $\PP n$ is said to be ``${\cal
O}_{\PP n}(d)$-numerically settled" if $deg A=h^{0}({\cal O}_{\PP
n}(d))$; in this case, $h^{0}({\cal I}_{A}(d))=0$ if and only if
$h^{1}({\cal I}_{A}(d))=0$. The scheme $Y_{n,d}$ is ${\cal O}_{\PP
n}(d)$-numerically settled for all $n,d$.

\a {\bf Remark 2.1} Let $A$ be a 0-dimensional ${\cal O}_{\PP
n}(d)$-numerically settled subscheme of $\PP n$, and assume
$h^{0}({\cal I}_{A}(d))=0$. Let $B\subseteq A$ and $C\supseteq A$ be
0-dimensional subschemes of $\PP n$; then, $h^{0}({\cal
I}_{C}(d))=0$, and $h^{1}({\cal I}_{B}(d))=0$, or equivalently,
$h^{0}({\cal I}_{B}(d))=deg A-deg B$.
\medskip \noindent Hence if we prove $h^{0}({\cal I}_{Y_{n,d}}(d))=0$ then we know that
$h^{1}({\cal I}_{Y_{n,d}}(d))=0$, and

\medskip  $h^{0}({\cal I}_{X_{s,n}}(d))=0$ for all $s>s_{n,d}$,
\medskip  $h^{1}({\cal I}_{X_{s,n}}(d))=0$ for all $s\leq s_{n,d}$.

\medskip \noindent Moreover, if $h^{0}({\cal I}_{Y_{n,d}}(d))=0$ then also $h^{0}({\cal
I}_{D}(d))=0$, where $D$ denotes a generic union of $X_{s_{n,d}}$,
of $\lfloor {r_{n,d}\over 2}\rfloor$ 2-jets and of $r_{n,d}-2\lfloor
{r_{n,d}\over 2}\rfloor$ simple points. In fact, we have
$h^{0}({\cal I}_{X_{s_{n,d}}}(d))=\deg (R_{n,d})= r_{n,d}$ and we
conclude by Lemma 1.8.
\par The same conclusion (i.e. $h^{0}({\cal I}_{D}(d))=0$) holds in the weaker assumption
that $h^{1}({\cal I}_{X_{s_{n,d}}}(d))=0$,  since in this case
$h^0({\cal I}_{X_{s_{n,d}}}(d))={d+n\choose n}-deg (X_{s_{n,d}})=
r_{n,d}$ and we get $h^{0}({\cal I}_{D}(d))=0$ by Lemma 1.8.

\a {\bf Theorem 2.2} {\it Suppose that for all $n\geq 5$, we have
$h^{1}({\cal I}_{X_{s_{n,3}}}(3))=0$ and $h^{0}({\cal
I}_{X_{s_{n,3}+1}}(3))=0$; then $h^{0}({\cal
I}_{Y_{n,d}}(d))=h^{1}({\cal I}_{Y_{n,d}}(d))=0$, for all $d\geq 4$,
$n\geq 4$.}

\medskip \Proof
Let us consider a hyperplane $H\subset \PP n$; we want a scheme $Z$
with support on $H$, made of $(2,3,n)$-schemes, and an integer
vector ${\bf p}$, such that the ``differential trace" $Tr^{\bf
p}_H(Z)\subset H$ is ${\cal O}_{\PP {n-1}}(d)$-numerically settled.

Let us consider $n\geq 5$ first. Since $0\leq r_{n-1,d} \leq 2n-2$,
we write $r_{n-1,d}=n\delta + 2h+\epsilon$, with $0\leq \epsilon
\leq 1$, $0\leq \delta \leq 1$ and $2h+\epsilon \leq n$.

We denote by $Z$ the zero dimensional scheme union of
$s_{n-1,d}+h+\epsilon + \delta$ (hence $\delta =0$ if $0\leq
r_{n-1,d} \leq n$, while $\delta =1$ if $n+1\leq r_{n-1,d} \leq
2n-2$), $(2,3,n)$-schemes
$Z_1,\dots,Z_{s_{n-1,d}+h+\epsilon+\delta}$, where each $Z_i$ is
supported at $P_i$ with direction $L_i$, and:
\medskip -  the $P_i$'s are generic on $H$, $i=1, \ldots , s_{n-1,d}+h+\epsilon+\delta$;
\medskip -  $L_i \subset H$ for $i=1, \ldots , s_{n-1,d}+h$;
\medskip - if $(\epsilon ,\delta)\neq (0,0)$, the corresponding lines
$L_{s_{n-1,d}+h+1}$, $L_{s_{n-1,d}+h+2}$ have
 generic directions in $\PP n$ (hence not contained in $H$.

In case $n=4$, instead, we write $r_{3,d}=2h+\epsilon$, with $0\leq
\epsilon \leq 1$, and $Z$ is given as before. Notice that in this
case $0\leq h \leq 3$, and it can appear only one line
$L_{s_{3,d}+h+1}$, not contained in $H$.

\noindent We want to use the Horace differential Lemma 1.5,  where
the role of the schemes $H$ and $Z$ appearing in the statement of
the Lemma are played  by  our hyperplane $H$ and the scheme $Z$ just
defined, and with:
\smallskip \b $W=A_{s_{n-1,d}+h+\epsilon+1}\cup \cdots \cup A_{s_{n,d}}\cup
R_{n,d}$,
\smallskip \b $S=A_1 \cup \dots \cup A_{s_{n-1,d}+h+\epsilon+\delta}$,
\smallskip \b ${\bf p}=(\underbrace {0,\dots,0}_{s_{n-1,d}}, \underbrace {1,\dots,1}_h,
\underbrace {2}_{\epsilon},\underbrace {0}_{\delta})$.

\smallskip \b so that $Tr_HW=\emptyset$ and $Res_HW=W$, and $Y_{n,d}=W\cup S$.
\b Notice that this construction is possible, since
$s_{n-1,d}+h+2\leq s_{n,d}$ (and even more than that): see Appendix
A, A.1.

\b In order to simplify notations, we set:
$$ T^j_i:=Tr^j_H(Z_i), \quad R^j_i:=Res^j_H(Z_i), \quad j=0,1,2,\quad
i=1,\dots,s_{n-1,d}+h+\epsilon+\delta,$$
$$T:=Tr_HW\cup Tr^{\bf p}_H(Z)= T^0_1\cup \dots \cup T^0_{s_{n-1,d}}\cup
T^1_{s_{n-1,d}+1}\cup \dots \cup T^1_{s_{n-1,d}+h}\cup
T^2_{s_{n-1,d}+h+\epsilon}\cup T^0_{s_{n-1,d}+h+\epsilon+\delta},$$
$$R:=Res_HW \cup Res^{\bf p}_H(Z)= W\cup R^0_1\cup \dots \cup R^0_{s_{n-1,d}}\cup
R^1_{s_{n-1,d}+1}\cup \dots \cup R^1_{s_{n-1,d}+h}\cup
R^2_{s_{n-1,d}+h+\epsilon}\cup R^0_{s_{n-1,d}+h+\epsilon+\delta}.$$

\noindent Observe that, by Remark 1.7 :

\medskip \b $T^0_1, \ldots , T^0_{s_{n-1,d}}$ are $(2,3,n-1)$-points in $H\cong \PP
{n-1}$, and $R^0_1, \ldots , R^0_{s_{n-1,d}}$ are 2-jets in $H$; \b
$T^1_{s_{n-1,d}+1}, ... , T^1_{s_{n-1,d}+h}$ are 2-jets in $H$ and
$R^1_{s_{n-1,d}+1}, ... , R^1_{s_{n-1,d}+h}$ are $(2,3,n-1)$-points
in $H$; \b $T^2_{s_{n-1,d}+h+\epsilon}$ is, when appearing, a simple
point of $H$, and $R^2_{s_{n-1,d}+h+\epsilon+\delta}$ is a  2-fat
point of $H$ doubled in a direction transversal to $H$;\b
$T^0_{s_{n-1,d}+h+\epsilon+\delta}$ is, when appearing, a 2-fat
point on $H$, and $R^0_{s_{n-1,d}+h+\epsilon}$ is a  2-fat point in
$\PP n$ with support on $H$.

\b We will also make use of the scheme:
$$B:= W\cup R^1_{s_{n-1,d}+1}\cup \dots \cup R^1_{s_{n-1,d}+h}\cup
R^2_{s_{n-1,d}+h+\epsilon}.$$

\b Let us consider the following four statements:
$${\bf Prop}(n,d): h^{0}({\cal I}_{Y_{n,d}}(d))=0; \qquad \quad {\bf Reg}(n,d):
h^{1}({\cal I}_{X_{s,n}}(d))=0\; {\rm and}\; h^{0}({\cal
I}_{X_{s,n}+1}(d))=0,$$
$${\bf Degue}(n,d): h^0({\cal I}_{R}(d-1)) = 0;\quad \qquad {\bf {Dime}}(n,d): h^0({\cal
I}_{T,H}(d)) =0.$$

\noindent If {\bf{ Degue}}$(n,d)$ and {\bf{ Dime}}$(n,d)$ are true,
we know that {\bf Prop}$(n,d)$ is true too, by Proposition 1.5.
\bigskip

For the first values of $n,d$, we will need an ``ad hoc"
construction, which is given by the following:
\bigskip
\noindent{\bf Lemma 2.3} {\it Let $d=4$ and $n\in \{4,5,6\}$, then
{\bf Prop}$(n,d)$ holds.}
\medskip
\noindent {\it Proof of the Lemma.}
\par $Case$ $n=4$.  Here we use the construction of
$R$ and $T$ described above, hence we need to show that {\bf{
Degue}}$(4,4)$ and {\bf{ Dime}}$(4,4)$ hold. Since $s_{3,4}=5$, and
$r_{3,4}=0$, $T$ is made of five generic $(2,3,3)$-points in $H\cong
\PP 3$, so {\bf{ Dime}}$(4,4)$ holds (i.e. $h^0(\PP
3,I_{T,H}(4))=h^0(\PP 3,I_{X_{5,3}}(4)=0$), e.g. see \ref {CGG1}.

In order to prove {\bf{ Degue}}$(4,4)$ we want to apply Lemma 1.2,
with $R$ made of five 2-jets plus the scheme $B = W$; hence we need
to show that $h^0({\cal I}_B(3))\leq 10$, while $h^0({\cal
I}_{Res_H(B)}(2))=0$. Since here $s_{4,4}=7=r_{4,4}$, while
$r_{3,4}=0$, we have that $B = W = Res_H(B)$ and it is given by
$A_6$ and $A_7$, plus $R_{4,4}$. Hence we have $h^1({\cal
I}_B(3))=0$, since $B$ is contained in the scheme made of 3 generic
$(2,3,4)$-points (which is known to have maximal Hilbert function,
by \ref {CGG1} or \ref{B}); $h^1({\cal I}_B(3))=0$ is equivalent to
saying that $h^0({\cal I}_B(3))=2s_{3,4}=10$, as required. Moreover
$h^0({\cal I}_B(2))=0$, since there is one only form of degree two
passing through two generic $(2,3,4)$-points in $\PP 4$, given by
the hyperplane containing the two double lines, doubled. Since the
support of $R_{4,4}$ is generic, we get $h^0({\cal I}_B(2))=0$. So
we have that {\bf{ Degue}}$(4,4)$ holds, and {\bf Prop}$(4,4)$ holds
too.

$Case$ $n=5$.  Here we need to use a different construction. We have
$s_{5,4}=11$, $r_{5,4}=5$, $s_{4,4}=7=r_{4,4}$. We want to use the
Horace differential Lemma 1.5 with $Z=Z_1\cup \dots \cup Z_8\cup
R_{5,4}$, where $Z_1, \dots , Z_8$ are $(2,3,5)$ schemes supported
at generic points of $H$ with direction $L_1, \dots, L_8 \subset H$,
and we specialize $R_{5,4}$ so that $R_{5,4}\subset H$, contained in
a generic $(2,3,4)$-scheme of $H$; with $W=A_9\cup A_{10} \cup
A_{11}$, and with ${\bf p}=(\underbrace {0,\dots,0}_{7}, 1, 0)$. \b
Hence $T=Tr_HW\cup Tr^{\bf p}_H(Z)= T^0_1\cup  T^0_2 \cup \dots \cup
T^0_7\cup T^1_8 \cup  R_{5,4}$ and $R=Res_HW \cup Res^{\bf p}_H(Z)=
W\cup R^0_1\cup R^0_2 \cup \dots \cup  R^0_7\cup R^1_8$.

We have that the ideal sheaf of $T^0_1\cup  T^0_2 \cup \dots \cup
T^0_7 \cup R_{5,4}$ has $h^1=0$ and $h^0=2$ in degree 4, by using
the previous case and the fact that $R_{5,4}$ is contained in a
$(2,3,4)$-point, so $h^0({\cal I}_{T,H}(4) ) = 0$ by Lemma 1.8,
since $T^1_8$ is a 2-jet in $H\cong \PP 4$. We also have $h^0({\cal
I}_{R}(3) )=0$. In fact,  let us denote by $U$ the scheme
$U=R^1_{8}\cup W$. In order to apply Lemma 1.9 (the $R^0_{i}$'s are
2-jets) to get $h^0({\cal I}_{R}(3) ) = 0$, we need to show that
$h^0({\cal I}_{Res_{H}U}(2) ) = 0$ and $h^1({\cal I}_{U}(3) ) = 0$.
Since $U$ is included in the union of four $(2,3,5)$-points, which
impose independent conditions in degree three (e.g. see \ref
{CGG1}), $h^1({\cal I}_{U}(3) ) = 0$ follows. Moreover, $Res_{H}(U)$
is made by three $(2,3,5)$-points, and again $h^0({\cal
I}_{Res_{H}U}(2) ) = 0$ is known by \ref {CGG1}.

Now, $h^0({\cal I}_{T,H}(4) ) = 0 = h^0({\cal I}_{R}(3) )$ imply
{\bf Prop}$(5,4)$ by Lemma 1.5, and we are done.

$Case$ $n=6$.  Here we have $s_{6,4}=16$, $r_{6,4}=2$,
 while $s_{5,4}=11$, $r_{5,4}= 5$.
We want to use the Horace differential Lemma 1.5 with $Z=Z_1\cup
\dots \cup Z_{13}\cup R_{6,4}$, where $Z_1, \dots , Z_{13}$ are
$(2,3,6)$ schemes supported at generic points of $H$ with direction
$L_1, \dots, L_{12} \subset H$, while $L_{13}$ is not in $H$, and we
specialize $R_{6,4}\subset H$, as a generic 2-jet in $H$; with
$W=A_{14}\cup A_{15} \cup A_{16}$,  and with ${\bf p}=(\underbrace
{0,\dots,0}_{11}, 1,2, 0)$. \b Hence $T=Tr_HW\cup Tr^{\bf p}_H(Z)=
T^0_1\cup  T^0_2 \cup \dots \cup T^0_{11}\cup T^1_{12}\cup T^2_{13}
\cup  R_{6,4}$ and $R=Res_HW \cup Res^{\bf p}_H(Z)= W\cup R^0_1\cup
R^0_2 \cup \dots \cup R^0_{11}\cup R^1_{12}\cup R^2_{13}$.

We have that $h^0({\cal I}_{T,H}(4) ) = 0$ by applying Lemma 1.1 and
the previous case.

We also have $h^0({\cal I}_{R}(3) )=0$. In fact, let us denote by
$U$ the scheme $U=R^1_{12}\cup R^2_{13}\cup W$. In order to apply
Lemma 1.9 (the $R^0_{i}$'s are 2-jets) to get $h^0({\cal I}_{R}(3) )
= 0$, we need to show that $h^0({\cal I}_{Res_{H}U}(2) ) = 0$ and
$h^1({\cal I}_{U}(3) ) = 0$.

Since $U$ is included in the union of five $(2,3,6)$-points, which
impose independent conditions in degree three (e.g. see \ref
{CGG1}), $h^1({\cal I}_{U}(3) ) = 0$ follows. Moreover, $Res_{H}(U)$
is made by three $(2,3,6)$-points plus a 2-fat point inside $H \cong
\PP 5$. Since there is only one form of degree two passing through
three generic $(2,3,6)$-points in $\PP 6$, given by the hyperplane
containing the three double lines, doubled, we get $h^0({\cal
I}_{Res_{H}U}(2) ) = 0$.

Now, $h^0({\cal I}_{T,H}(4) ) = 0 = h^0({\cal I}_{R}(3) )$ imply
{\bf Prop}$(6,4)$ by Lemma 1.5, and we are done. \hfill \qed

\a Now we come back to the proof of the Theorem for the remaining
values of $n,d$; we will work  by induction on both $n,d$ in order
to prove statement {\bf Prop}$(n,d)$ for $n\geq 4$, $d\geq 5$ and
for $n\geq7$, $d=4$. We divide the proof in 7 steps.

\medskip
\noindent {\it Step 1.} The induction is as follows: we suppose that
{\bf Prop}$(\nu,\delta)$ is known for all $(\nu,\delta)$ such that
$4\leq \nu < n$ and $4\leq \delta \leq d$ or $4\leq \nu \leq n$ and
$4\leq \delta < d$ and  we prove that {\bf Prop}$(n,d)$ holds.

The initial cases for the induction are given by Lemma 2.2, and we
will also make use of the fact that ${\bf Reg}(n,3)$ with $n\geq 4$
and ${\bf Reg}(3,d)$ with $d\geq 4$ hold respectively by assumption
and by [B], while, by [CGG], we know everything about the Hilbert
function of generic $(2,3,n)$-schemes when $d=2$.

We will be done if we prove that {\bf{ Degue}}$(n,d)$ and {\bf{
Dime}}$(n,d)$ hold for $n\geq 4$, $d\geq 5$ and for $n\geq7$, $d=4$.

\medskip
\noindent {\it Step 2.} Let us prove {\bf{ Dime}}$(n,d)$. Notice
that $T$ is ${\cal O}_{\PP {n-1}}(d)$-numerically settled in $H\cong
\PP {n-1}$, hence {\bf{ Dime}}$(n,d)$ is equivalent to $h^1({\cal
I}_{T,H}(d)) =0$.

The scheme $T$ is the generic union of $X_{s_{n-1,d}}$ with $h$
2-jets, of $\epsilon$ simple points and of $\delta$ 2-fat points,
where $2h+\epsilon+n\delta =r_{n-1,d}$. Then {\bf {Dime}}$(n,d)$
holds for $n\geq 5$ and $d\geq 4$ since we are assuming that {\bf
Prop}$(n-1,d)$ is true and the union of $h$ 2-jets, $\epsilon$
simple points and of $\delta$ 2-fat points can specialize to
$R_{n-1,d}$ (see Lemma 1.10).

For $n=4$ and $d\geq 5$, {\bf {Dime}}$(4,d)$ holds, since we know
that $h^1({\cal I}_{X_{s_{3,d}}}(d))=0$ by [B] and in this case $T$
is the generic union of $X_{s_{3,d}}$ with $h$ 2-jets and $\epsilon$
simple points so we can apply Lemma 1.8.

\medskip \noindent {\it Step 3.} We are now going to prove
{\bf{ Degue}}$(n,d)$. Since the scheme $R$ is the union of the
scheme $B$ and of $s_{n-1,d}\;$ 2-jets lying on $H$ (see definitions
of $R$ and $B$ above), we can use Lemma 1.9 $ii)$. Hence, in order
to prove that $\dim(I_{R} )_{d-1} = 0$, i.e. that
{\bf{Degue}}$(n,d)$ holds, it is enough to prove that
$(I_{Res_{H}(B)})_{d-2}=0$ and that $\dim (I_{B})_{d-1}\leq
2s_{n-1,d}$.

\medskip
\noindent {\it Step 4.} Let us show that$(I_{Res_{H}(B)})_{d-2}=0$.
We set $t_{n,d}:=s_{n,d} - s_{n-1,d}-h-\epsilon-\delta$.
 The scheme $Res_{H}(B)$ is given by $W$ plus, if $\epsilon =1$, one
2-fat point contained in $H$, plus, if $\delta =1$, one simple point
in $H$. $W$ is the generic union of $R_{n,d}$ with $t_{n,d}$
$(2,3,n)$-points. Let $I$ denote the ideal of these $t_{n,d}$
$(2,3,n)$-points; if we show that $I_{d-2}=0$, then also
$(I_{Res_{H}(B)})_{d-2}=0$.

The idea is to prove that our $(2,3,n)$-points are ``too many" to
have $I_{d-2}\neq 0$ since they are more than $s_{n,d-2}+1$; the
only problem with this procedure is that there are cases (when
$d-2=2$ or $3$) where $I_{d-2}$ may not have the expected dimension,
so those cases have to be treated in advance.

First let $d=4$ (and $n\geq 7$); if we show that $t_{n,4}>{n\over
2}$, then we are done, since $(I_{X_{s,n}})_{2}=0$ for $s> {n\over
2}$, by \ref {CGG}, Prop 3.3.  The inequality $t_{n,4} > {n\over 2}$
is treated in Appendix A, A.2, and proved for $n\geq 7$, as
required.

Now let $d=5$ and $n=4$; here we have that $s_{4,3}+1=4$, but
actually there is one cubic hypersuface through four
$(2,3,4)$-points in $\PP 4$; nevertheless, since $t_{4,5}= 14-8-0-0
= 6$, and it is known (see \ref {CGG}or \ref B) that
$(I_{X_{6,4}})_{3}=0$, we are done also in this case.

Eventually, for $d= 5$, $n\geq 5$, or in the general case $d\geq 6$,
$n\geq 4$, if we show that $t_{n,d}\geq s_{n,d-2} +1$, the problem
reduces to the fact that $(I_{X_{s_{n,d-2}}+1})_{d-2}=0$. If $d=5$,
we know that $(I_{X_{s_{n,3}}+1})_3=0$ by hypothesis, while for
$d\geq 6$ we can suppose that $(I_{X_{s_{n,d-2}}+1})_{d-2}=0$ by
induction on $d$.

The inequality $t_{n,d}\geq s_{n,d-2} +1$ is discussed in Appendix
A, A.1, and proved for all the required values of $n,d$.

Thus the condition $(I_{Res_{H}(B)})_{d-2}=0$ holds.

\medskip
\noindent {\it Step 5.} Now we have to check that $\dim
(I_{B})_{d-1}\leq 2s_{n-1,d}$. Since $\deg Y_{n,d}= h^0({\cal
O}_{\PP n}(d))$ and $\deg T= h^0({\cal O}_{\PP {n-1}}(d))$, then
$\deg R= h^0({\cal O}_{\PP {n}}(d-1))$. The scheme $R$ is the union
of the scheme $B$ and of $s_{n-1,d}\;$ 2-jets lying on $H$, so $\deg
R= \deg B+2s_{n-1,d}$. Hence $\dim (I_{B})_{d-1}\leq 2s_{n-1,d}$ is
equivalent to $h^1({\cal I}_B(d-1))=0$ (and to $\dim (I_{B})_{d-1}=
2s_{n-1,d}$).

Let us consider the case $n\geq 5$ first. Let $Q$ be the scheme
$Q=Z_{s_{n-1,d}+1}\cup \dots \cup Z_{s_{n,d}+h+\epsilon+\delta}\cup
A_{s_{n,d}+h+\epsilon+\delta+1}\cup \dots \cup A_{s_{n,d}}\cup
A_{s_{n,d}+1}$, where $A_{s_{n,d}+1}$ is a $(2,3,n)$ scheme
containing  $R_{n,d}$. We have that $B$ is contained in the scheme
 $Q$, which is composed by $s_{n,d}-s_{n-1,d}+1$ generic
$(2,3,n)$-points (notice that $2h+\epsilon+\delta \leq n+1$, so
$Z_{s_{n-1,d}+1},\dots , Z_{s_{n,d}+h+\epsilon+\delta}$ are generic,
since only the first $h$ of the lines $L_i$ are in $H$).

\par The generic union of $s_{n,d-1}$ generic
$(2,3,n)$-points in $\PP n$ is the scheme $X_{s_{n,d-1}}$; by
induction, or by hypothesis if $d-1=3$, we have $h^{1}({\cal
I}_{X_{s_{n,d-1}}}(d-1))=0$. Since $s_{n,d}-s_{n-1,d}+1\leq
s_{n,d-1}$ (see {\it Step 6}), then $B\subset Q \subset
X_{s_{n,d-1}}$ and we conclude by Remark 2.1 that $h^{1}({\cal
I}_{B}(d-1))=0$.

\medskip
\noindent {\it Step 6.} We now prove the inequality:
$s_{n,d}-s_{n-1,d}+1\leq s_{n,d-1}$ ($n\geq 5$). \b We have $\deg Q
= \deg B + 2h + \epsilon + n\delta + (2n+1-r_{n,d}) $, in fact in
order to ``go from $B$ to $Q$", we have to add  a 2-jet to each of
the $R^1_{i}$ ($h$ in number), a simple point to
$R^{2}_{s_{n-1,d}+h+\epsilon}$ if $\epsilon =1$, a 2-fat point of
$H$ if $\delta = 1$ and something of degree $(2n+1-r_{n,d})$ to
$R_{n,d}$.

Since $r_{n,d}\geq 0$ and $2h+\epsilon+n\delta = r_{n-1,d}\leq
2n-2$, we have: $\deg Q=(2n+1)(s_{n,d}-s_{n-1,d}+2 ) \leq \deg
(B)+2n-2+2n+1=\deg (B)+4n-1$.

Notice that $\deg(Y_{n,d-1})=\deg (B) + 2s_{n-1,d}$, so we have:
$(2n+1)(s_{n,d}-s_{n-1,d}+1 )\leq \deg(Y_{n,d-1}) -
2s_{n-1,d}+4n-1$.

If we prove that $4n-1-2s_{n-1,d}\leq 0$, we obtain:
$(2n+1)(s_{n,d}-s_{n-1,d}+1 )\leq \deg(Y_{n,d-1}) =
(2n+1)s_{n,d-1}$, and we are done.

The computations to get  $4n-1-2s_{n-1,d}\leq 0$ can be found in
Appendix A.3.

\medskip
\noindent {\it Step 7.} We are only left to prove that $h^1({\cal
I}_B(d-1))=0$ in case $n=4$ ( $d\geq 5$).

Recall that now $r_{3,d} = 2h +\epsilon \leq 6$, with $0\leq h\leq
3$, $0\leq \epsilon \leq 1$.  If $r_{3,d}\leq 4$, we can apply the
same procedure as in step 5, since the part of the scheme $Q$ with
support on $H$ is generic in $\PP 4$.  Hence we only have to deal
with $r_{3,d} = 5,6$.

The case $r_{3,d} = 5$ does not actually present itself; this can be
checked by considering that
$$ {d+3 \choose 3} = {(d+3)(d+2)(d+1) \over 6} = 7s_{3,d} +
r_{3,d} \Rightarrow (d+3)(d+2)(d+1) = 42s_{3,d} + 6r_{3,d}.$$ Hence
if $r_{3,d}=5$, we get $42s_{3,d} + 30 = 7(6s_{3,d} + 4) + 2$, but
it is easy to check that $(d+3)(d+2)(d+1)$  never gives a remainder
of 2, modulo 7.

Thus we are only left with the case $r_{3,d} = 6$, when $h=3$ and
$\epsilon =0$.  In this case we have $d \equiv  3$ (mod 7), hence
$d\geq 10$; it is also easy to check that $r_{3,d-1} = 3$ in this
case.

We can add $2s_{3,d}$ generic simple points to $B$, in order to get
a scheme $B'$ which is ${\cal O}_{\PP 4}(d-1)$-numerically settled,
so now $h^1({\cal I}_B(d-1))=0$ is equivalent to $h^0({\cal
I}_{B'}(d-1))=0$ (by Remark 2.1).

We want to apply Horace differential Lemma again in order to prove
$h^0({\cal I}_{B'}(d-1))=0$; so we will  define appropriate schemes
$Z_B$, $W_B$ and an integer vector ${\bf q}$, such that conditions
$a)$ and $b)$ of Proposition 1.5 apply to them, yielding $h^0({\cal
I}_{B'}(d-1))=0$.

Consider the scheme $Z_B\subset \PP 4$, given by  $s_{3,d-1}-1$
$(2,3,4)$-schemes in $\PP 4$, such that their support is at generic
points of $H$, and only for the last one of them the line $L_i$ is
not in $H$.  Let $W_B\subset \PP 4$ be given by $2s_{3,d}$ generic
simple points, $s_{4,d} - s_{3,d} - s_{3,d-1}-2$ generic
$(2,3,4)$-schemes, three generic $(2,3,3)$-schemes in $H\cong \PP
3$, and the scheme $R_{4,d}$. Let also ${\bf q}=(\underbrace
{0,\dots,0}_{s_{3,d-1}-3},  \underbrace {1}_{1},\underbrace
{2}_{1})$.

Let $T_B = Tr_H(W_B)\cup Tr_H^{\bf p}(Z_B)= X_{s_{3,d-1}}\cup E \cup
F$,   and  $R_B = Res_H(W_B)\cup Res_H^{\bf q}(Z_B)$.

We have that $E$ and $F$ are, respectively, a 2-jet and a simple
point in $H$ (they give the ``remainder scheme" of degree 3, to get
that $T_B$ is ${\cal O}_{\PP 3}(d-1)$-numerically settled).

The scheme $R_B$ is the union of $2s_{3,d}$ generic simple points,
$s_{4,d} - s_{3,d} - s_{3,d-1}-2$ generic $(2,3,4)$-schemes, the
scheme $R_{4,d}$, $s_{3,d-1}$ 2-jets in $H$, a $(2,3,3)$-scheme in
$H$ and  a  2-fat point of $H$ doubled in a direction transversal to
$H$.

If we show that $h^0({\cal I}_{R_B}(d-2))=0 = h^0({\cal
I}_{T_B,H}(d-1))$, then we are done by Proposition 1.5.

We have $h^0({\cal I}_{T_B,H}(d-1))=0$, since $T_B$ is ${\cal
O}_{\PP 3}(d-1)$-numerically settled, and is given by the union of
$X_{s_{3,d-1}}$ (whose ideal sheaf has $h^1=0$ in degree $d-1$ by
[B]) with a 2-jet and a simple point, so we can apply Lemma 1.8.

In order to show that $h^0({\cal I}_{R_B}(d-2))=0$ we want to
proceed as in Step 5, i.e  by applying Lemma 1.9, since $R_B$,
 is made of $s_{3,d-1}-3$ 2-jets union the $2s_{3,d}$ generic simple points
 and a scheme that we denote by $R'_B$. We will be done if we show that $h^0({\cal
 I}_{Res_H(R_B)}(d-3))=0$ and $h^1({\cal I}_{R'_B}(d-2))=0$.

 The first condition will follow if $s_{4,d}-s_{3,d}-s_{3,d-1}-2\geq s_{4,d-3}$, the
second condition (since $R'_B$ is contained in
 the union of $s_{4,d}-s_{3,d}-s_{3,d-1}+1$ generic $(2,3,4)$-schemes) if
$s_{4,d}-s_{3,d}-s_{3,d-1}+1\leq
 s_{4,d-2}$.

 Both inequalities are proved in Appendix, A.4.

\hfill \qed

\bigskip
Thanks to some ``brute force" computation by {\bf COCOA}, we are
able to prove:
\bigskip {\bf Corollary 2.4} {\it For $4\leq n \leq
9$, we have:
\par i) $h^{1}({\cal I}_{X_{s_{n,3}}}(3))=0$ and $h^{0}({\cal I}_{X_{s_{n,3}+1}}(3))=0$,
except for $n=4$, in which case we have $h^{0}({\cal
I}_{X_{s,4}}(3))=0$ for $s\geq 5$.
\par ii) $h^{0}({\cal I}_{Y_{n,d}}(d))=h^{1}({\cal I}_{Y_{n,d}}(d))=0$, for $d\geq 4$.}

\medskip \Proof \ Part i) comes from direct computations using CoCoA ([{\bf
CO}]). Note that $s_{4,3}=3$ and that $h^{0}({\cal
I}_{X_{4,4}}(3))=h^{1}({\cal I}_{X_{4,4}}(3))=1$, see [{\bf CGG1}].
\par \noindent  Part ii) comes by applying the Theorem and part i).

\bigskip
Coming back to the language of secant varieties, Theorem 2.2 and
Corollary 2.4 give:
\bigskip {\bf Corollary 2.5} {\it  If Conjecture 1 is true for $d=3$, then it is true for
all $d\geq 4$. Moreover, for $n \leq 9$, Conjecture 1 holds.}

\bigskip
\bigskip

\a {\bf 3. On Conjecture 2a. The case $n=2$.}

\bigskip
In this section we prove Conjecture 2a for $n=2$.

We want to use the fact that $\sigma_s (O_{k,n,d})$ is defective if
at a generic point its tangent space does not have the expected
dimension; actually (see \ref {BCGI}) this is equivalent to the fact
that for generic $L_i\in R_1$, $F_i\in R_k$,
$R=\kappa[x_0,...,x_n]$, $i=1,...,s$ the vector space
$<L_1^{d-k}R_k,L_1^{d-k-1}F_1R_1,...,L_s^{d-k}R_k,L_s^{d-k-1}F_sR_1>$
does not have the expected dimension.

Via inverse systems this reduces to the study of $(I_Y)_d$, where $Y
= Z_1\cup ...\cup Z_s$ is a certain 0-dimensional scheme in $\PP n$.
Namely, the scheme $Y$ is supported at $s$ generic points
$P_1,...,P_s \in \PP n$, at each of them $\deg(Z_i) = {k+n \choose
n}+n$, and $I_{P_i}^{k+2}\subset I_{Z_i}\subset I_{P_i}^{k+1}$ (see
Lemma 1.2).

\bigskip

When working in $\PP 2$, we can specialize the $F_i$'s to be of the
form $\Pi_i^k$, where $\Pi_i$ is a generic linear form through
$P_i$.  In this way we get a scheme $\overline{Y}=
\overline{Z}_1\cup ...\cup \overline{Z}_s$, and the structure of
each $\overline{Z}_i$ is $((k+2)P_i\cap L_i^2)\cup (k+1)P_i$, where
the line $L_i$ is ``orthogonal" to $\Pi_i=0$, i.e. if we put
$P_i=(1,0,0)$, $\Pi_i=x_1$  and $L_i= \{x_2=0\}$, the ideal is of
the form: $((x_1,x_2)^{k+2}+(x_2)^2)\cap (x_1,x_2)^{k+1} =
(x_1^{k+2}, x_1^{k+1}x_2, x_1^{k-1}x_2^2,...,x_2^{k+1})$.

Notice that the forms in $I_{\overline Z_i}$ have multiplicity at
least $k+1$ at $P_i$ and they meet $L_i$ with multiplicity at least
$k+2$; moreover the generic form in $I_{\overline Z_i}$ has $L_i$ at
least as a double component of its tangent cone at $P_i$.

\medskip

 When $F\in I_{\overline{Z}_i}$ and we speak of its ``tangent cone" at
$P_i$, we mean (with the choice of coordinates above) either the
form in $\kappa[x_1,x_2]$ obtained by putting $x_0=1$ in $F$ and
considering the (homogeneous) part of minimum degree thus obtained,
or also the scheme (in $\PP 2$) defined by such a form.

When we will say that $L_i$ is a ``simple tangent" for $F$, we will
mean that $L_i$ is a reduced component for the tangent cone to $F$
at $P_i$.

\bigskip
The strategy we adopt to prove Conjecture 2a is the following: if
$(I_Y)_d$ does not have the expected dimension, i.e. $h^0({\cal
I}_Y(d))h^1({\cal I}_Y(d))\neq 0$, then the same happens for ${\cal
I}_{\overline Y}(d)$; hence Conjecture 2a would be proved if we show
that whenever $\dim (I_{\overline Y})_d$ is more than expected, then
 $h^1({\cal I}_ {X}(d))> max \{0, \deg(Y) - {d+n \choose n} \}\;$  or
  $\;h^0({\cal I}_ {T}(d))>max \{0, {d+n \choose n}-\deg(Y)
\}$, where
$$X:= (k+1)P_1\cup ...\cup (k+1)P_s \subset \PP 2; \ \ \
  T:= (k+2)P_1\cup \ldots \cup (k+2)P_s \subset \PP 2 . $$

 \bigskip
The following easy technical Bertini-type lemma and its corollary will be of use
in the sequel.
\bigskip

\noindent {\bf Lemma 3.1}  {\it Let $F$, $G $ be linearly independent polynomials in $\kappa[x]$. Then for almost any $a\in
\kappa$, $F+aG$ has at least one simple root.}

{\it Proof}. Let $M$ be the greatest common divisor of $F$ and $G$
with $F=MP$, $G=MQ$. Let us consider $PQ'-QP'$, where $P'$ and
$Q'$ are the derivatives of $P$ and $Q$, respectively . Since $P$
and $Q$ have no common roots, it easily follows that
 $PQ'-QP'$  cannot be identically zero.
 \par

 For any
$\beta \in \kappa$ which is neither a root for  $PQ'-QP'$, nor for $M$, nor for $Q$, let

 $$a=a(\beta) := - {P(\beta)\over Q(\beta) },$$
 so
 $(F+aG)(\beta)=M(\beta)(P+aQ)(\beta)= 0,$ and
 $(F+aG)'(\beta)= ( M'(P+aQ) + M (P' +aQ'))  (\beta)=(M(P' +aQ'))
(\beta) =(M(P' - {P(\beta)\over Q(\beta) }Q')) (\beta)=({M \over
Q})(\beta)  (QP'-PQ')(\beta) \neq 0 ,$
 hence $\beta$ is a simple root for $F+aG$.
Since $\beta$ assumes almost every value in $\kappa$, so does $a(\beta)$. \hfill \qed

\bigskip

{\bf Corollary 3.2} {\it Let $P=(1,0,0)\in \PP 2$. Let $f,g\in
(I_P^{k+1})_d$, and  $f,g \notin (I_P^{k+2})_d$.  Assume that  $f,g$, have different  tangent cones at $P$. Then for almost any $a\in
\kappa$,  $f+ag$ has
at least one simple tangent at $P$.}

\medskip

{\it Proof}. The Corollary follows immediately from Lemma 3.1 by
de-homogenising the tangent cones to $f,g$ at $P$ to get two non-zero and non-proportional polynomials $F,G\in \kappa[x]$. \hfill \qed

\bigskip

It  will be handy to introduce the following definitions.
\par
\bigskip

{\bf Definition 3.3 } Let $P\in \PP 2$ and $L$ be a line $L$ through
$P$. We say that a scheme supported at one point is of type ${Z'}$
if its structure is $(k+1)P \cup ((k+2)P \cap L)$, and that it is of
type $\overline{Z}$ if its structure is $(k+1)P \cup ((k+2)P \cap
L^2)$.
 \par
 We will say that a union of schemes of types ${Z'}$ and/or
 $\overline{Z}$ is generic if the points of their support and the
 relative lines are generic.

\bigskip

 The following lemma is the key to prove
Conjecture 2a:

\bigskip

{\bf Lemma 3.4} {\it Let $\overline{Y}= \overline{Z}_1\cup ...\cup
\overline{Z}_s \subset \PP 2$ be a union of $s$ generic schemes of
type $\overline{Z}$, then either:

{\rm (i)} $(I_{\overline Y} )_d = (I_T)_d$;

or

{\rm (ii)} $\dim (I_{\overline Y} )_d = \dim (I_X)_d - 2s$.}

\medskip

{\it Proof}.
 Notice that by  the genericity of the points and of the lines, the Hilbert function of a scheme
 with support on $P_1, \ldots, P_s $, formed
by $t$  schemes of  type  $\overline{Z}$, by $t'$  schemes of  type
${Z'}$ and by $s-t-t' $ fat points of multiplicity $(k+1)$ depends
only on $s, t $ and $t'$.
\par

Let  $W_t$ be a scheme  formed by $t$  schemes of  type
$\overline{Z}$  and by $s-t $ fat points of multiplicity $(k+1)$.
Let
$$ \tau = \max \{ t \in \Bbb N | \dim (I_{W_t })_d = \dim (I_X)_d - 2t  \} .
$$
\medskip

If $\tau = s$, we have $W_{s} = \overline Y$ and  $\dim (I_{W_s})_d
= \dim (I_X)_d - 2s$ , hence (ii) holds.
\par Let  $\tau  <s$: we will prove that $(I_{\overline Y} )_d = (I_T)_d$.
Let $W$ be the scheme
$$ W = W_\tau = \overline{Z}_1\cup \ldots \cup \overline{Z}_{  \tau  } \cup
 (k+1)P_{  \tau  +1 }\cup \ldots \cup (k+1)P_s .
$$
and let
$$ W'_{(j)} = \overline{Z}_1\cup \ldots \cup \overline{Z}_{  \tau  } \cup
(k+1)P_{  \tau  +1 }\cup \ldots \cup {Z'}_{  j } \cup \ldots  \ldots
\cup (k+1)P_s , \ \ \ \  \ \tau  +1 \leq j \leq s  ,$$

$$ W''_{(j)} = \overline{Z}_1\cup \ldots \cup \overline{Z}_{  \tau  } \cup
(k+1)P_{  \tau  +1 }\cup \ldots \cup  \overline{Z}_{  j } \cup
\ldots  \ldots \cup (k+1)P_s , \ \ \ \  \ \tau  +1 \leq j \leq s
,$$ that is  $W'_{(j)}$, respectively  $W''_{(j)}$, is the scheme
obtained from  $W$ by substituting the fat point  $(k+1)P_{ j }$
with a scheme of type $Z'$, respectively $\overline Z$, so
$$ W \subset W'_{(j)} \subset W''_{(j)} ,
$$
and  $\deg  W'_{(j)}  = \deg  W +1$,  $\deg  W''_{(j)}  = \deg  W
+2$ (for $ \tau = s-1$, $ W''_{(s)}= \overline Y $).
\par

If  $ (I_{W''_{(j)}  })_d =0$, then trivially $(I_{\overline Y} )_d
=(I_T)_d=0$ and we are done. So assume that $ (I_{W''_{(j)}  })_d
\neq 0$.
\par

By the definition of $\tau$ we have that $ \dim (I_{W''_{(j)} })_d
>   \dim (I_X)_d - 2(\tau +1) = \dim (I_{W })_d -2$, hence we get
$$ 0 \leq \dim (I_{W  })_d - \dim (I_{W''_{(j)}   })_d  \leq 1 .
$$

Let us consider the two possible cases. \par

\medskip

 {\it Case 1:  }
 $\dim (I_{W  })_d - \dim (I_{W'_{(j)}  })_d =0 \  , \  \   \tau  +1 \leq j \leq s  $.\par

In this case we have   $(I_{W  })_d =  (I_{W'_{(j)}  })_d$ . This
means that every form $F \in  (I_{W })_d  $ meets the line $L_{j}$
with multiplicity at least $k+2$; but since the line $L_{j}$ is
generic through $P_{j}$, this yields that every line through $P_{j}$
is met with multiplicity at least $k+2$, hence
$$(I_{W})_d \subset
(I_{P_{j}}^{k+2})_d , \ \ {\rm for } \ \  \tau +1 \leq j \leq s .
\eqno (1) $$
\par

In particular, we have that
$$ (I_{W})_d = (I_{W''_{(s)}})_d .  \eqno (2)
$$

Now consider the schemes
$$ W_{(i,s)}=
\overline {Z}_1 \cup  \ldots \cup \overline {Z}_{i-1} \cup (k+1)P_i \cup \overline {Z}_{i+1}  \cup \ldots \cup
\overline{Z}_{  \tau  } \cup
 (k+1)P_{  \tau  +1 } \cup \ldots \cup (k+1)P_{s-1}  \cup \overline Z_s  \ ,
\ \ \ \  \ 1 \leq i \leq \tau ,
$$
$$ W'_{(i,s)}=
\overline {Z}_1 \cup \ldots \cup \overline {Z}_{i-1} \cup Z'_i \cup \overline {Z}_{i+1} \cup \ldots \cup
\overline{Z}_{  \tau  } \cup
 (k+1)P_{  \tau  +1 } \cup \ldots \cup (k+1)P_{s-1}  \cup \overline Z_s  \ ,
\ \ \ \  \ 1 \leq i \leq \tau ,
$$
i.e.  $W_{(i,s)}$  is the scheme obtained from  $W$ by substituting
the fat point $(k+1)P_i $ to the scheme $\overline Z_i $  and a
scheme $ \overline Z _s$, of type $ \overline Z $, to the fat point
$(k+1)P_{s }$,  while $W'_{(i,s)}$  is the scheme obtained from
$W_{(i,s)}$  by substituting a scheme $Z'_i$, of type $  Z' $, to
  the fat point $(k+1)P_i $.

  The schemes  $W_{(i,s)}$ and $W$ are made of  $\tau$ schemes of type $\overline Z$ and $s-\tau$ $(k+1)$-fat points;
the schemes $W'_{(i,s)}$ and $W'_{(s)}$ are made of  $\tau$ schemes
of type $\overline Z$, $s-\tau-1$ $(k+1)$-fat points and one scheme
of type $Z'$. This yields that:
 $$\dim ( I_{  W_{(i,s)} } )_d = \dim ( I_ W )_d=  \dim ( I_{  W'_{(s)} } )_d = \dim ( I_{  W'_{(i,s)} } )_d .
 $$
 Hence every form  $F \in  (I_{ W_{(i,s)}  })_d  $ meets the generic  line $L_{i}$ with
multiplicity at least $k+2$, thus we get
$$(I_{ W_{(i,s)} })_d \subset
(I_{P_{i}}^{k+2})_d , \ \ {\rm for } \ \  1 \leq i \leq \tau . \eqno
(3) $$ and from this and (2) we have
$$ (I_{W_{(i,s)} })_d = (I_{W''_{(s)}})_d = (I_{W})_d .  \eqno (4)
$$
By (1), (3) and (4) it follows that $(I_{W})_d  = (I_{T})_d$, hence,
since $W \subset \overline Y \subset T$, we get (i).
\par

 \bigskip

 {\it Case 2} :  $\dim (I_{W  })_d - \dim (I_{W'_{(j)}  })_d =1 \ , \ \ \   \tau  +1 \leq j \leq s $.\par
 In this case we have
 $$  \dim (I_{W'_{(j)}  })_d = \dim (I_{W''_{(j)}  })_d .
 $$
 Let $F \in  (I_{W'_{(j)}  })_d = (I_{W''_{(j)}  })_d $; hence $L_j$ appears with multiplicity two in the tangent cone of $F$.
If  $F \not\in  (I^{k+2}_{P_{j}  })_d $, then let $L'_j$ be a
generic line not in the tangent cone of $F$ at $P_{j}$. By
substituting the line  $L'_j$ to $L_j$ in the construction of
$W'_{(j)}$, we get another form
 $G \in (I  _{W})_d $, $G \not\in  (I^{k+2}_{P_{j}  })_d $, with the double line $L'_j$ in its tangent cone.
 Then, by Corollary 3.2, the generic form $F+aG$ has a simple tangent at $P_j$, and this is a contradiction
 since a generic choice of the line $L_j$ should yield $(I_{W'_{(j)}  })_d =  (I_{W''_{(j)}  })_d$. Hence $F \in  (I^{k+2}_{P_{j}  })_d $,
 for $ \tau+1 \leq j \leq s$.

With an argument like the one we used  in {\it Case 1},  we also get
that $F \in  (I^{k+2}_{P_{j}  })_d $ for $1 \leq j \leq \tau$, and
(i) easily follows. \hfill \qed

\bigskip

 \par
Now we are ready to prove Conjecture 2a.
\bigskip
{\bf Theorem 3.5} {\it The secant variety $\sigma_s (O_{k,2,d})$ is
defective if and  only if one of the following holds:
\par
(i) $h^1({\cal I}_ {X}(d))> max \{0, \deg(Y) - {d+n \choose n}
\}\;$, or
\par
(ii) $\;h^0({\cal I}_ {T}(d))>max \{0, {d+n \choose n}-\deg(Y) \}$.}

\medskip

{\it Proof}. Since if $Y$ is defective in degree $d$, then
$\overline Y$ is, but, by Lemma 3.4, either
  $\dim (I_{\overline Y})_d = \dim (I_X)_d - 2s$, hence
   $$h^1({\cal I}_ {X}(d)) =
h^1({\cal I}_ {\overline Y}(d))  -2s
   > max \{0, \deg( {\overline Y}) - {d+n \choose n} \} = max \{0, \deg(Y) - {d+n \choose n} \}  ,$$
  or  $(I_{\overline Y})_d =
(I_T)_d$, hence
 $$ h^0({\cal I}_ {T}(d)) =  h^0({\cal I}_ {\overline Y}(d)) >
 max \{0, {d+n \choose n}-\deg(\overline Y) \}=
 max \{0, {d+n \choose n}-\deg( Y)\}.$$

\hfill \qed

\bigskip

\bigskip\bigskip\centerline {{\bf APPENDIX: Calculations}}

\bigskip

\noindent {\bf A.1} \quad We want to prove that (for $n\geq 4$ and
$d\geq 6$ or for $n\geq 5$ and $d= 5$):

$$ s_{n,d} -  s_{n-1,d} - h - \epsilon -\delta -1 \geq  s_{n,d-2}$$

Recall:

$ s_{n,d}(2n+1) +  r_{n,d} =  {n+d\choose d}$;\quad  $
s_{n-1,d}(2n-1) + r_{n-1,d} =  {n+d-1\choose d}$;\quad  $
s_{n,d-2}(2n+1) + r_{n,d-2} = {n+d-2\choose d-2}$.
\medskip
Hence our inequality becomes:
\medskip
$ {1\over 2n+1}\left[{n+d\choose d} -  r_{n,d}\right] - {1\over
2n-1}\left[{n+d-1\choose d} -  r_{n-1,d}\right] - h - \epsilon
-\delta -1 - {1\over 2n+1}\left[{n+d-2\choose d-2} -
r_{n,d-2}\right] \geq 0$
\medskip
By using binomial equalities and reordering this is:
\medskip
\noindent $ {1\over 2n+1}\left[{n+d-1\choose d}+{n+d-2\choose d-1} +
{n+d-2\choose d-2}\right] - {1\over 2n-1}{n+d-1\choose d}
+{r_{n-1,d}\over 2n-1} - h -\epsilon -\delta -1 - {1\over
2n+1}{n+d-2\choose d-2} +{1\over 2n+1}( r_{n,d-2}- r_{n,d}) \geq 0$

\medskip
i.e.
\medskip
  ${1\over 2n+1}{n+d-2\choose d-1} - {2\over (2n+1)(2n-1)}{n+d-1\choose
  d}+{r_{n-1,d}\over 2n-1} - h -\epsilon -\delta -1 +{1\over 2n+1}(
r_{n,d-2}- r_{n,d})\geq 0$
\medskip
By using binomial equalities again:
\medskip
  ${1\over 2n+1}{n+d-2\choose d-1} - {2\over (2n+1)(2n-1)}\left[{n+d-2\choose
  d}+{n+d-2\choose d-1}\right]+{r_{n-1,d}\over 2n-1} - h -\epsilon
  - \delta -1 +{1\over 2n+1}( r_{n,d-2}- r_{n,d})\geq 0$
\medskip
i.e.
\medskip
${1\over 2n+1}{n+d-2\choose d-1}(1- {2\over 2n-1})- {2\over
(2n+1)(2n-1)}{n+d-2\choose
  d}+{r_{n-1,d}\over 2n-1} - h -\epsilon -\delta -1 +{1\over 2n+1}(
r_{n,d-2}- r_{n,d})\geq 0$
\medskip
i.e.
\medskip
${n+d-2\choose d-1}{[2n(d-1)-3d+2]\over d(4n^2-1)}+{r_{n-1,d}\over
2n-1} - h -\epsilon -\delta -1 +{1\over 2n+1}( r_{n,d-2}-
r_{n,d})\geq 0$
\medskip
Now, ${r_{n-1,d}\over 2n-1}\geq 0$, while  $h + \epsilon +\delta
\leq {n\over 2}$, and $r_{n,d-2}- r_{n,d}\geq -2n$, i.e. ${1\over
2n+1}( r_{n,d-2}- r_{n,d})\geq -{2n\over 2n+1}\geq -1$, so our
inequality holds if:

\medskip ${n+d-2\choose d-1}{[2n(d-1)-3d+2]\over
d(4n^2-1)}-{n\over 2}-2\geq 0$

  It is quite immediate to check that the right hand side is an increasing
function in $d$, e.g. by writing it as follows:
\medskip
 ${n+d-2\choose n-1}[2n-3-{2n+2\over d}] -
({n\over 2}+2)(4n^2-1)\geq 0.$
\medskip
i.e.
\medskip
 ${n+d-2\choose n-1}[2n-3-{2n+2\over d}] - 2n^3-8n^2+{n\over 2}+2\geq 0.$

\bigskip
Let us consider the case $d=6$ first; our inequality becomes:
\medskip
${n+4\choose 5}{(10n-16)\over 6} - 2n^3-8n^2+{n\over 2}+2\geq 0.$
\medskip
i.e.
\medskip
${(n+4)(n+3)(n+2)(n+1)n(5n-8)\over 360}- 2n^3-8n^2+{n\over 2}+2\geq
0.$
\medskip
i.e.
\medskip
${(n+4)(n+3)(n+2)(n+1)n(5n-8)-20n^2(n+2)\over 360}+{n\over 2}+2\geq
0.$
\medskip
i.e.
\medskip
${n(n+2)\over 360}[(n+4)(n+3)(n+1)(5n-8)-720n]+{n\over 2}+2\geq 0.$
\medskip
Which, for $n\geq 4$, is easily checked to be true. Hence we are
done for $n\geq 4$, $d\geq 6$.

\bigskip
Now let us consider the case $d=5$; our inequality becomes:
\medskip
${n+3\choose 4}{(8n-13)\over 5} - 2n^3-8n^2+{n\over 2}+2\geq 0.$
\medskip
i.e.
\medskip
${(n+3)(n+2)(n+1)n(8n-13)\over 120}- 2n^3-8n^2+{n\over 2}+2\geq 0.$
\medskip
i.e.
\medskip
$(n^4+6n^3+11n^2+6n)(8n-13)- 240n^3-960n^2+60n+240\geq 0.$
\medskip
i.e.
\medskip
$8n^5 + 35n^4 - 230n^3 - 1015n^2 - 18n + 240\geq 0.$
\medskip
i.e.
\medskip
$n^3(8n^2+35n- 230-{1015\over n}-{18\over n^2}+{240\over n^3})\geq
0.$
\medskip
Which, for $n\geq 6$, holds. So we are left to prove our inequality
for $d=5=n$; in this case we have:  $s_{5,5}= [{272\over 11}] = 24$,
$s_{4,5}= [{126\over 9}] = 14$ and $r_{4,5}= 0$, hence $h= \epsilon
= 0$, while $s_{5,3}= [{56\over 11}] = 5$;  so:\quad $ s_{5,5} -
s_{4,5} -1 \geq s_{5,3}$ becomes: $24 -14 -1 \geq 5$, which holds.

\hfill \qed

\bigskip

\a {\bf A.2} We want to prove that, for all $n\geq 7$:

$$s_{n,4}-s_{n-1,4}-h-\epsilon -\delta > {n\over 2}$$

i.e.

 ${n+4 \choose 4}/(2n+1)-r_{n,4}/(2n+1)-{n-1+4\choose
4}/(2n-1)+r_{n-1,4}/(2n-1)-h-\epsilon - \delta >{n\over 2}$

 i.e.

${(n+4)(n+3)(n+2)(n+1)\over 24(2n+1)}-{(n+3)(n+2)(n+1)n\over
24(2n-1)}-{n\over 2}-{r_{n,4}\over (2n+1)} + {r_{n-1,4}\over (2n-1)}
-h-\epsilon - \delta
>0$

Now:

 ${r_{n,4}\over (2n+1)}\leq {2n\over (2n+1)}<1$, hence $-{r_{n,4}\over (2n+1)}>-1$;

 $r_{n-1,4}\geq 0$;

 and $h+\epsilon +\delta \leq {n\over 2}$, i.e. $-h-\epsilon -\delta \geq - {n\over 2}$.

 Therefore we get:

${(n+3)(n+2)(n+1)\over 24} \cdot [{(n+4)\over (2n+1)}- {n\over
(2n-1)}]-{n\over 2}-{r_{n,4}\over (2n+1)}+ {r_{n-1,4}\over
(2n-1)}-h-\epsilon -\delta >$

$ {(n+3)(n+2)(n+1)\over 24} \cdot [{n+4\over 2n+1}-{n\over
2n-1}]-{n\over 2}-{n\over 2}-1=$

$={(n+3)(n+2)(n+1)\over 24} \cdot {[(2n-1)(n+4)-n(2n+1)]\over
(2n+1)(2n-1)}- n-1=$

$(n+1)\left[{(n+3)(n+2)(3n-2)\over 12(4n^2-1)}-1\right] > 0$

i.e.

$(n+3)(n+2)(3n-2)-12(4n^2-1)>0$

i.e.

$3n^3-35n^2+8n>0$

which is true for $n\geq 12$.

Let us check the cases $n=7,8,9,10,11$.

If $n=7$ we have: $s_{7,4}=[{1\over 15}{11\choose 4}]=22$ (with
$r_{7,4}=0$); $s_{6,4}=16$, since $ {10 \choose 4}=210=16\cdot
13+2$, so $r_{6,4}=2$ and $h=1$, $\epsilon = \delta =0$.

Our inequality becomes: $22-16-1> 7/2$, which holds.

If $n=8$ we have: $s_{8,4}=[{1\over 15}{12\choose 4}]=33$ (with
$r_{8,4}=0$); $s_{7,4}=22$,  $r_{7,4}=0$ and $h=\epsilon = \delta
=0$.

Our inequality becomes: $33-22> 4$, which holds.

If $n=9$ we have: $s_{9,4}=[{1\over 15}{13\choose 4}]=47$ (with
$r_{9,4}=10$); $s_{8,4}=33$, and  $h=\epsilon = \delta =0$.

Our inequality becomes: $47-33> 9/2$, which holds.

If $n=10$ we have: $s_{10,4}=[{1\over 15}{14\choose 4}]=66$ (with
$r_{10,4}=11$) ; $s_{9,4}=47$, and  $h=5$, $\epsilon = \delta =0$.

Our inequality becomes: $66-47-5> 5$, which holds.

If $n=11$ we have: $s_{10,4}=[{1\over 15}{15\choose 4}]=91$;
$s_{10,4}=66$, and  $h=5$, $\epsilon = 1$, $\delta =0$.

Our inequality becomes: $91-66-5-1> 11/2$, which holds. \hfill \qed

\bigskip

\noindent {\bf A.3} We want to prove that, for $d\geq 5$, $n\geq 4$
or $d=4$, $n\geq 7$:
 $$4n-1\leq 2s_{n-1,d}.    \qquad \qquad \qquad (*) $$
\medskip
Since $r_{n-1,d}\leq 2n-2$, it is enough to prove that:
\medskip
 ${2\over 2n-1}\left[ {n-1+d\choose n-1}-2n+2
\right]\geq 4n-1$
 which is:
\medskip
 $ {n-1+d\choose n-1}
\geq {(4n-1)(2n-1) \over 2}+2n-2$
 that is:
\medskip
 $ {n-1+d\choose n-1}
\geq 4n^2-n-{3\over 2}   \qquad  \qquad \qquad (**)$

 \medskip which is surely true if
 \medskip $ {n-1+d\choose n-1}
\geq 4n^2-n$ is true.
\medskip
Notice that the function  $ {n-1+d\choose n-1}$ is an increasing
function in $d$. For $d=4$, the inequality becomes:
\medskip ${n(n^3+6n^2+11n+6)\over 24}\geq 4n^2-n$, which can be written:

\medskip $n^3+6n^2+11n+6\geq 96n-24$, i.e.

\medskip $n^3+6n^2-85n+30\geq 0$ which is surely true if the following is true:

\medskip $n^2+6n-85\geq 0$.
The last one is verified for $n\geq 8$, so we are done for $d=4$ and
$n\geq 8$.

\a If  $(n,d)=(7,4)$, $s_{n-1,d}=16$ since $ {10 \choose
4}=210=16\cdot 13+2$, and $(*)$ becomes: $4\cdot7-1\leq 2\cdot 16$
which is true.

Since the function  $ {n-1+d\choose n-1}$ is an increasing function
in $d$, we have proved the initial inequality for $d\geq4$ and
$n\geq 8$.

For $d=5$ $(**)$ becomes: $n^5 + 10n^4 + 35n^3 - 430n^2 + 144n + 120
\geq 0$ which is true for $ n= 5,6,7$.  We have hence proved the
initial inequality for $d\geq 5$ and $n\geq 5$.

\a If  $(n,d)=(4,5)$, $s_{n-1,d}=8$ since $ {8 \choose 3}=8\cdot 7$,
and $(*)$ becomes: $4\cdot4-1\leq 2\cdot 8$ which is true.

For $d=6$  $(**)$ becomes: $n(n+1)(n+2)(n+3)(n+4)(n+5)-120(6)(4n^2
-n-1)\geq 0$ which is true for $ n= 4$.  We conclude that the
initial inequality is true for $d\geq 5$ and $n\geq 4$.

\bigskip

{\noindent {\bf A.4} We want to show that (for $d\geq 10$):
$s_{4,d}-s_{3,d}-s_{3,d-1}-2 \geq s_{4,d-3}$ and
$s_{4,d}-s_{3,d}-s_{3,d-1}+1 \leq s_{4,d-2}$  .

The first inequality is equivalent to:

$$\left[ {1\over 9}{d+4\choose 4}\right] - {1\over 7}{d+3\choose 3} +
{6\over 7} - {1\over 7}{d+2\choose 3} + {3\over 7} -2 \geq \left[
{1\over 9}{d+1\choose 4}\right]$$

which follows if:
\medskip

$ {1\over 9}{d+4\choose 4} -{1\over 9}{d+1\choose 4}  \geq {1\over
7}{d+3\choose 3}+{1\over 7}{d+2\choose 3}- {9\over 7}+4$

\medskip
i.e.
\medskip
$ {d+1\over 9}{[(d+4)(d+3)(d+2)-d(d-1)(d-2)]\over 24}  \geq {1\over
7}\left({(d+1)(d+2)(2d+3)\over 6}\right)+ {19\over 7}$
\medskip
i.e.
\medskip
$ {d+1\over 9}{(12d^2+24d+24)\over 24}  \geq {1\over
42}(d+1)(d+2)(2d+3) + {19\over 7}$
\medskip
i.e.
\medskip
$ {(d^2+2d+2)\over 3}  \geq {2d^2+7d+6\over 7} + {114\over 7(d+1)}$
\medskip
i.e.
\medskip
$ d^2-7d-4  \geq  {342\over d+1}$

\medskip
Which is easily checked to hold for $d\geq 10$.

\bigskip
Now let us consider the second inequality, which is equivalent to:

$$\left[ {1\over 9}{d+4\choose 4}\right] - {1\over 7}{d+3\choose 3} +
{6\over 7} - {1\over 7}{d+2\choose 3} + {3\over 7} +1 \leq \left[
{1\over 9}{d+2\choose 4}\right]$$

which follows if:
\medskip

$ {1\over 9}{d+4\choose 4} -{1\over 9}{d+2\choose 4}  \leq {1\over
7}{d+3\choose 3}+{1\over 7}{d+2\choose 3}- {9\over 7}-3$
\medskip
i.e.
\medskip
$ {(d+1)(d+2)\over 9}{[(d+4)(d+3)-d(d-1)]\over 24}  \leq {1\over
7}\left({(d+1)(d+2)(2d+3)\over 6}\right)-{30\over 7}$
\medskip
i.e.
\medskip
$ {(d+1)(d+2)\over 9}{(8d+12)\over 24}  \geq {1\over
42}(d+1)(d+2)(2d+3) - {30\over 7}$
\medskip
i.e.
\medskip
$ {1\over 9}  \geq {1\over 7} - {180\over 7(d+1)(d+2)(2d+3)}$

\medskip
Which is easily checked to hold for $d\geq 10$.

\bigskip\bigskip\centerline {{\bf REFERENCES}}

\medskip\noindent[{AH1}]: J. Alexander, A. Hirschowitz. {\it Polynomial
interpolation in several variables.} J. of Alg. Geom. {\bf 4}
(1995), 201-222.

\medskip\noindent[{AH2}]: J. Alexander, A. Hirschowitz. {\it An asymptotic vanishing theorem for generic
unions of multiple points.} Invent. Math.  {\bf 140} (2000), 303-325.

\medskip\noindent [{B}]: E. Ballico, {\it On the
secant varieties to the tangent developable of a Veronese variety},
  J. Algebra  {\bf 288} (2005), 279--286.

\medskip\noindent[{BBCF}] E. Ballico, C.Bocci, E.Carlini,
C.Fontanari.  {\it Osculating spaces to secant varieties.}  Rend.
Circ. Mat. Palermo {\bf 53 } (2004), 429-436.

\medskip\noindent [{BF}]: E. Ballico, C.Fontanari,{\it On the
secant varieties to the osculating variety of a Veronese surface},
Central Europ. J. of Math. {\bf 1} (2003), 315-326.

\medskip\noindent [{BC}]: A.Bernardi, M.V.Catalisano.
{\it Some defective secant varieties to osculating varieties of Veronese surfaces.}
Collect. Math. {\bf 57}, n.1, (2006), 43-68.

\medskip\noindent[{BCGI} ]: A.Bernardi, M.V.Catalisano, A.Gimigliano, M.Id\`a.
{\it  Osculating Varieties of Veronesean and their higher secant
varieties. }  Canadian J. of  Math. {\bf  59}, (2007), 488-502.

\medskip\noindent[{BO}]: M.C.Brambilla, G.Ottaviani, {\it On the
Alexander-Hirschowitz theorem.} J. Pure and Appl. Algebra {\bf 211}
(2008), 1229-1251.

\medskip\noindent [{CGG}]: M.V.Catalisano, A.V.Geramita, A.Gimigliano.
{\it On the Secant Varieties to the Tangential Varieties of a
Veronesean.}  Proc. A.M.S.  {\bf 130} (2001), 975-985.

\medskip\noindent [{Ch1}]: K.Chandler. {\it A brief proof of a
maximal rank theorem for generic double points in projective space.}
Trans. Amer. Math. Soc. {\bf 353} (2001), 1907-1920.

\medskip\noindent [{Ch2}]: K.Chandler. {\it Linear systems of cubics singular at general points of projective space.}
Composition Math. {\bf 134} (2002), 269-282.

\medskip\noindent [{Ge}]: A.V.Geramita. {\it Inverse Systems of Fat
Points}, Queen's Papers in Pure and Applied Math.  {\bf 102}, {\it
The Curves Seminar at Queens', vol. X} (1998).

\medskip\noindent [{IK}]: A.Iarrobino, V.Kanev. {\it Power Sums,
Gorenstein algebras, and determinantal loci.} Lecture Notes in Math.
{\bf 1721}, Springer, Berlin, (1999).

\medskip\noindent [{I}]: A.Iarrobino. {\it Inverse systems of a symbolic
algebra III: Thin algebras and fat points.} Compos. Math. {\bf 108}
(1997), 319-356.

\medskip\noindent [{J}]: J.P.Jouanoulou. {\it Th\'eoremes de Bertini
et applications}. Progress in Math. {\bf 42}, Birkauser,  Boston,
MA, (1983).

\medskip\noindent [{Se}]: B. Segre, {\it Un'estensione delle variet\`a di
Veronese ed un principio di dualit\`a per le forme algebriche I and
II}. Rend. Acc. Naz. Lincei (8) {\bf 1} (1946), 313-318 and 559-563.

\medskip\noindent [{Te}]: A.Terracini. {\it Sulle} $V_k$ {\it per cui la
variet\`a degli} $S_h$ $(h+1)${\it -seganti ha dimensione minore
dell'ordinario.} Rend. Circ. Mat. Palermo {\bf 31} (1911), 392-396.
\medskip
\noindent [{ W}]: K.Wakeford. {\it On canonical forms.} Proc. London
Math. Soc. (2) {\bf 18} (1919/20). 403-410.
\par
\medskip

\bigskip \bigskip

\b {\it A.Bernardi, Dip. Matematica, Univ. di Bologna, Italy, email:
abernardi@dm.unibo.it}
 \medskip
\b {\it M.V.Catalisano, Dip. Matematica, Univ. di Genova, Italy,
e-mail: catalisano@diptem.unige.it}
\medskip
\b {\it A.Gimigliano, Dip. di Matematica and C.I.R.A.M., Univ. di
Bologna, Italy, e-mail: gimiglia@dm.unibo.it}
\medskip
\b {\it M.Id\`a, Dip. di Matematica, Univ. di Bologna, Italy,
e-mail: ida@dm.unibo.it}

\end